# HOMEOMORPHISMS OF 3-MANIFOLDS AND
# THE REALIZATION OF NIELSEN NUMBER

Boju Jiang[1], Shicheng Wang[1,2], and Ying-Qing Wu[2]

Abstract. The Nielsen Conjecture for Homeomorphisms asserts that any homeomorphism $f$ of a closed manifold is isotopic to a map realizing the Nielsen number of $f$, which is a lower bound for the number of fixed points among all maps homotopic to $f$. The main theorem of this paper proves this conjecture for all orientation preserving maps on geometric or Haken 3-manifolds. It will also be shown that on many manifolds all maps are isotopic to fixed point free maps.

The proof is based on the understanding of homeomorphisms on 2-orbifolds and 3-manifolds. Thurston's classification of surface homeomorphisms will be generalized to 2-dimensional orbifolds, which is used to study fiber preserving maps of Seifert fiber spaces. Maps on most Seifert fiber spaces are indeed isotopic to fiber preserving maps, with the exception of four manifolds and orientation reversing maps on lens spaces or $S^3$. It will also be determined exactly which manifolds have a unique Seifert fibration up to isotopy. These informations will be used to deform a map to certain standard map on each piece of the JSJ decomposition, as well as on the neighborhood of the decomposition tori, which will make it possible to shrink each fixed point class to a single point, and remove inessential fixed point classes.

## Table of Contents

1991 *Mathematics Subject Classification.* 57N10, 55M20, 57M50.
[1] Partially supported by an NSFC grant
[2] Research at MSRI supported in part by NSF grant #DMS 9022140.

Typeset by $\mathcal{AMS}$-TEX

## §0. Introduction

Consider a map $f$ on a compact polyhedron $M$. Denote by $Fix(f)$ the set of fixed points of $f$. Two points $x, y \in Fix(f)$ are equivalent (or Nielsen equivalent) if there is a path $\gamma$ from $x$ to $y$ such that $\gamma$ and $f \circ \gamma$ are rel $\partial\gamma$ homotopic. It can be shown that this is an equivalence relation, so it divides $Fix(f)$ into equivalence classes, called the Nielsen classes. For each fixed point class one can define an index. A fixed point class is essential if its index is non-zero. The Nielsen number $N(f)$ of $f$ is then defined as the number of essential fixed point classes in $Fix(f)$. The important property of $N(f)$ is that it is a homotopy invariant of maps. In particular it gives a lower bound of the number of fixed points for all maps homotopic to $f$. One is referred to [J1] for an introduction to the Nielsen fixed point theory.

It is known that if $M$ is a compact manifold of dimension $\geq 3$, then $N(f)$ can be realized as the number of fixed points for some map homotopic to $f$, see [Ki, We, Br]. This result is not true for continuous maps on hyperbolic surfaces [J2]. For homeomorphisms of manifolds, Nielsen conjectured that if $f$ is a surface homeomorphism then the Nielsen number is realized by the fixed point number of some homeomorphism isotopic to $f$. A proof of this conjecture was announced by Jiang [J3] and Ivanov [Iv], and given in detail by Jiang and Guo in [JG]. The following conjecture of Nielsen remains an open problem for manifolds of dimension 3 and 4.

**The Nielsen Conjecture for Homeomorphisms.** *If $M$ is a closed manifold, and $f : M \to M$ is a homeomorphism, then $f$ is isotopic to a homeomorphism $g$ whose fixed point number equals $N(f)$.*

This is trivial for 1-manifolds, and the work of Jiang and Guo [JG] proved it for 2-manifolds, which uses the Thurston's classification of surface homeomorphisms. Kelly [Ke] gave a proof of this conjecture for homeomorphisms on manifolds of dimension at least 5, using the fact that any null-homotopic simple closed curve in such manifolds bounds an embedded disc.

The main result of this paper is to prove this conjecture for orientation preserving homeomorphisms on closed orientable 3-dimensional manifolds which are either Haken or geometric. Thurston's Geometrisation Conjecture asserts that these include all irreducible, orientable, closed 3-manifolds.

**Theorem 9.1.** *Suppose $M$ is a closed orientable manifold which is either Haken or geometric, and $f : M \to M$ is an orientation preserving homeomorphism. Then* *f can be isotoped to a homeomorphism g with $\#Fix(g) = N(f)$.*

There are eight geometries for 3-manifolds. A 3-manifold in the theorem can be decomposed along a canonical set of tori into pieces, each of which admits one of these geometries. This decomposition is called the Jaco-Shalen-Johannson torus decomposition, or simply the JSJ decomposition. We will isotop the homeomorphism on each individual piece $M_i$ to a "standard map", which has the property that no two components of the fixed point set of $f|_{M_i}$ are equivalent on that piece, and the behavior of the map near each fixed point component are well understood. Some major difficulties arise when we try to put the pieces together, because two fixed points being equivalent is a global property: Two components of the fixed point set of $f|_{M_i}$ may be inequivalent in $M_i$, but yet they may be equivalent in the whole manifold. Therefore the homeomorphisms $f|_{M_i}$ should not only satisfy the above conditions, but also some extra conditions to guarantee that after gluing the pieces together we still get a homeomorphism with good properties. This justifies the sophisticated definition of standard maps given in Section 5. It can then be shown that if the restrictions of the homeomorphism on the pieces are standard maps in that sense, then the global homeomorphism indeed has the property that each fixed point class is connected. The local behavior of the homeomorphism near each fixed point component will then enable us to shrink that component to a single point, and delete it if it is an inessential class.

For many manifolds, a much stronger result holds: All homeomorphisms on those manifolds are isotopic to fixed point free maps. Define an orbifold $X(M)$ to be *small* if it is a sphere with a total of at most three holes or cone points, or a projective plane with a total of 2 holes or cone points; otherwise it is *big*.

**Theorem 9.2.** *Let $M$ be a closed orientable 3-manifold. Then any orientation preserving homeomorphism $f$ on $M$ is isotopic to a fixed point free homeomorphism, unless some component of the JSJ decomposition of $M$ is a Seifert fiber space with big orbifold.*

In particular, the theorem applies to the following manifolds: (1) all Seifert fiber space whose orbifold is a sphere with at most 3 cone points or a projective plane with at most two cone points. This covers all manifolds admitting a $S^3$ or $S^2 \times S^1$ geometry, and four of the six Euclidean manifolds. (2) All manifolds obtained by gluing hyperbolic manifolds together along boundary tori. This includes all closed hyperbolic manifolds. (3) All Sol manifolds.

Our proof strongly depends on the knowledge of homeomorphisms on geometric manifolds, as well as homeomorphisms on 2-dimensional orbifolds. In Section 1 we will generalize Thurston's classification theorem of hyperbolic surface homeomorphisms to homeomorphisms of hyperbolic orbifolds. We define periodic, reducible and pseudo-Anosov maps on hyperbolic orbifolds $X$ in a natural way, in the sense that if $F$ is a surface covering of $X$, and $\widetilde{f}$ is a lifting of a map $f$ on $X$, then $f$ is pseudo-Anosov if and only if $\widetilde{f}$ is.

**Theorem 1.4.** *Suppose $X$ is a hyperbolic orbifold. Then any homeomorphism $f : X \to X$ is isotopic to either a periodic map, a reducible map, or a pseudo-Anosov map.*

A similar result is proved in Section 2 for homeomorphisms on Euclidean orbifolds. These results are crucial to our proof because the orbifold of a Seifert fiber space $M$ has a natural orbifold structure, and any fiber preserving homeomorphism on $M$ induces a homeomorphism of the orbifold.

In Section 3 we study the problem of which manifold has a Seifert fibration such that every orientation preserving homeomorphism is isotopic to fiber preserving one. This property is known to be true for all Seifert fiber space which is not covered by $S^2 \times S^1$, $S^3$ or $T^3$ [Sc, Theorem 3.9], as well as for those manifolds which has a sphere orbifold with cone points of order $(2, 3, p)$, $p \geq 5$ or $(3, 3, k)$, $k \geq 2$ [BO, Proposition 3.1]. In Section 3 we will show that this is true for all but four Seifert fiber spaces.

**Theorem 3.11.** *Suppose $M$ is a compact orientable Seifert fiber space which is not $T^3$, $M_{P(2,2)}$, $S^1 \times D^2$, or $T \times I$. Then there is a Seifert fibration $p : M \to X(M)$, so that any orientation preserving homeomorphism on $M$ is isotopic to a fiber preserving homeomorphism with respect to this fibration. Moreover, if $M$ is not a lens space or $S^3$, then the result is true for all homeomorphisms on $M$.*

A strong version of the converse is true: Not only that none of the four manifolds listed has a "universal" invariant Seifert fibration, but also there exist homeomorphisms on these manifolds that do not preserve any Seifert fibration. As a corollary, we also completely determine all Seifert fiber spaces which have a unique Seifert fibration *up to isotopy*, see Corollary 3.12 for more details.

Two of the exceptional manifolds in Theorem 3.11 are closed. We need to find certain representatives for homeomorphisms on these manifolds separately. The manifold $T^3$ is easy: Each homeomorphism on $T^3$ is isotopic to a linear map. The manifold $M_{P(2,2)}$ is a Euclidean manifold with orbifold a projective plane with two cones of order 2. This manifold will be treated in Section 4. It is shown that any homeomorphism on it is isotopic to an isometry with respect to certain Euclidean metric.

Section 5 gives the definition of some basic concept about fixed point theory, including Nielsen class, its index, and the Nielsen number of a map. We will also define two types of standard maps on 3-manifolds, and prove some useful lemmas.

The remaining part of the paper is to use the results in the earlier sections to give a proof of Theorem 9.1. Given a generic 3-manifold, the Jaco-Shalen-Johannson decomposition cuts $M$ into hyperbolic and Seifert fibered pieces. Using Theorem 1.4, one can further cut the Seifert fibered components into pieces so that on each piece $M_i$ the homeomorphism $f$ induces an orbifold homeomorphism which is either periodic or pseudo-Anosov. The main result of Section 6 is to show that if $f$ is such a homeomorphism on a piece $M_i$, then it can be isotopic to a standard map.

In Section 7, we will study torus fiber preserving homeomorphisms on torus bundles over 1-orbifolds. This includes twisted $I$-bundles over the Klein bottle, the union of two such manifolds along their boundary, $T \times I$, and torus bundle over $S^1$. An important lemma is that if $f$ is already standard on the boundary of a $T \times I$, then it can be isotoped rel $\partial(T \times I)$ to a standard map. This is useful in putting standard maps on the pieces together to form a global standard map. Also, the study of homeomorphisms on torus bundles over $S^1$ proves Theorem 9.1 for Sol manifolds. The theorem is proved in Section 8 for five other types of geometric manifolds, one is the hyperbolic manifolds, all the other four types are Seifert manifolds with spherical or Euclidean orbifolds.

The final proof of Theorem 9.1 is given in Section 9. With the previous results we can isotop the homeomorphism so that it is standard on each piece, including the collars of the decomposing tori. It will then be shown that in this case the global map will have the property that each fixed point class is connected. The behavior of the homeomorphism in a neighborhood of each fixed point component is well understood, which enables us to shrink each fixed point class to a single point, and remove it if it is inessential. This will complete the proof of the main theorem. Theorem 9.2 is also proved in this section.

We refer the reader to [Sc] for basic definitions and properties of the eight 3-manifold geometries, which is crucial in understanding this paper. Basic concepts about 3-manifold topology can be found in [Ja, He], in particular the Jaco-Shalen-Johannson decomposition of Haken manifolds and Seifert fiber spaces. See [Th1, Th2] for discussions on hyperbolic 3-manifolds.

We would like to thank Peter Scott for some helpful conversations about maps on Seifert fiber spaces.

**Notations and Conventions.** Given a set $X$ in $M$, we use $N(X)$ to denote a regular neighborhood of $X$, use $|X|$ to denote the number of components in $X$, and use $\#X$ to denote the number of points in $X$.

If $f : M \to M$ is a map, we use $Fix(f)$ to denote the set of fixed points of $f$, and use $N(f)$ to denote the Nielsen number of $f$.

A map $f$ on a Seifert fiber space $M$ is a *fiber preserving map* if it maps each fiber to another fiber. Identifying each fiber of $M$ to a point, we get a set $X(M)$, which has a natural 2-dimensional orbifold structure [Sc]. If $f$ is fiber preserving, it induces a map $\widehat{f} : X(M) \to X(M)$. An isotopy $h_t$ ($t \in I$) of a Seifert fiber space is *fiber preserving* if each $h_t$ is a fiber preserving map. It is a *fiber-wise isotopy* if $h_t$ maps each fiber to itself.

Unless otherwise stated, all maps on 3-manifolds and 2-orbifolds are assumed to be homeomorphisms.

§1. Classification of maps on hyperbolic orbifolds

We refer the reader to [Sc] for definitions about 2-dimensional orbifolds, their Euler characteristics, and other basic concepts and properties. All orbifolds in this section are assumed compact. Recall that an orbifold is good if it is covered by surface. Scott [Sc, Theorem 2.3] gave a list of all bad orbifolds without boundary. It is in fact the list of all bad orbifolds: for any orbifold $X$ with nonempty boundary, its double $\widehat{X}$ along $\partial X$ is clearly not in that list, thus covered by some surface, hence so is $X$ itself. An orbifold is hyperbolic if it has negative Euler Characteristic. A hyperbolic orbifold is covered by a hyperbolic surface, and admits hyperbolic structure with totally geodesic boundary. A reflector circle or reflector arc will be called a reflector. Denote by $S(X)$ the set of singular points of $X$. Thus each component of $S(X)$ is either a cone point, a reflector, or a union of reflector arcs joined at reflector corners.

**Definition 1.1.** A *homeomorphism* of an orbifold $X$ is a map $f : X \to X$ such that (1) $f$ is a homeomorphism of the underlying topological space, and (2) $f$ preserves the orbifold structure, i.e., it maps a cone point (resp. reflector corner) to a cone point (resp. reflector corner) of the same angle, and it maps a reflector to a reflector.

Similarly, an *isotopy* on $X$ is a continuous family of orbifold homeomorphisms $f_t$, $0 \leq t \leq 1$. In this case $f_t$ is called an isotopy from $f_0$ to $f_1$.

Unless otherwise stated, all orbifold maps below are orbifold homeomorphisms.

**Definition 1.2.** A set of mutually disjoint curves $C = \cup C_i$ in the interior of a hyperbolic orbifold $X$ is an *admissible set* if

(1) Each $C_i$ is either a circle disjoint from $S(X)$, or is a reflector circle, or is an arc with interior disjoint from $S(X)$ and with ends on cone points of order 2.

(2) Each component of $X - \text{Int} N(C)$ is a hyperbolic orbifold.

It is easy to see that if $C$ is an admissible set on $X$, then its lifting $\widetilde{C}$ on a compact covering surface $\widetilde{X}$ is an admissible set of $\widetilde{f}$ in the usual sense [FLP], i.e., it consists of mutually disjoint, mutually non-parallel, essential simple closed curves on $\widetilde{X}$.

**Definition 1.3.** Suppose $f : X \to X$ is a map on a hyperbolic orbifold.

(1) $f$ is *periodic* if $f^n = id$ for some $n$.

(2) $f$ is *reducible* if $f(C) = C$ for some admissible set $C$. The set $C$ is called a *reducing set* of $f$.

(3) $f$ is a *pseudo-Anosov* map if $S(X)$ has only cone points, and $f : X' \to X'$ is a pseudo-Anosov map in the sense of Thurston, where $X'$ is the punctured surface $X - S(X)$.

Figure 1.1

A pseudo-Anosov map has two transverse measured foliations, $(\mathcal{F}_h, \tau_h)$ and $(\mathcal{F}_v, \tau_v)$. A measured foliation on the orbifold $X$ is also a measured foliation of the underlying surface, except that a cone point of $X$ may be a singular point of the foliations with prong number 1, see Figure 1.1.

Thurston's classification of surface maps [Th3] says that any map on a hyperbolic surface is isotopic to either a periodic map, a reducible map, or a pseudo-Anosov map. For a proof of this theorem on orientable surfaces, see [FLP, Be, CB]. When the surfaces are non orientable, this is proved in [Wu]. The following theorem generalizes this theorem to maps on hyperbolic orbifolds.

**Theorem 1.4.** *Suppose $X$ is a hyperbolic orbifold. Then any homeomorphism $f : X \to X$ is isotopic to either a periodic map, a reducible map, or a pseudo-Anosov map.*

*Proof.* If $X$ has reflector circles, let $C$ be the union of such. Then $f(C) = C$, and $X - \mathrm{Int} N(C)$ has the same Euler characteristic as that of $X$, so $f$ is a reducible map.

If $X$ has reflector arcs but no reflector circles, let $F$ be the underlying surface of $X$, and let $C' = C'_1 \cup \ldots \cup C'_n$ be the components of $\partial F$ containing some reflector arcs. Clearly we have $f(C') = C'$. Isotop $C'$ into the interior of $F$, we get a set of curves $C = C_1 \cup \ldots \cup C_n$ in $X$. By an isotopy we may assume that $f(C) = C$. The set $C$ cuts $X$ into $X_1 \cup \ldots \cup X_n \cup Y$, where $X_i$ is an annulus between $C'_i$ and $C_i$, and $Y$ is the remaining component. Since each $X_i$ is an annulus containing some reflector arcs, it is a hyperbolic orbifold. If $Y$ is also hyperbolic, then $C$ is an admissible set and we are done. So assume that $Y$ has nonnegative Euler characteristic. There are several cases.

CASE 1: $Y$ is a Möbius band or an annulus with $\partial Y = C$. Replace $C$ by the central curve of $Y$.

CASE 2: $Y$ is a disk with two cone points of order 2 and no other cone points. Replace $C = \partial Y$ by an arc in Y connecting the two cone points.

CASE 3: $Y$ is either an annulus with one boundary on $\partial X$, or a disk with at most one cone point. Then $f$ can be isotoped to a periodic map.

In all 3 cases $f$ can be isotoped to a periodic map or a reducible map. Therefore the theorem holds if $X$ has some reflector arcs.

Now assume that $X$ has no reflectors. Let $X' = X - S(X)$, and let $f'$ be the restriction of $f$ to $X'$. If $f'$ is isotopic to a periodic or pseudo-Anosov map, then so is $f$, and we are done. So assume $f'$ is isotopic to a reducible map $g'$ with reducing curves $C'$. Let $g$ be the corresponding map on $X$. Then $g$ is isotopic to $f$ via an orbifold isotopy. If each component of $X - \mathrm{Int}N(C')$ is a hyperbolic orbifold, then $g$ is a reducible map and we are done. If some component $X_i$ of $X - \mathrm{Int}N(C')$ is not hyperbolic, it must be a disk $D$ with two cone points of order 2. Replace $\partial D$ with an arc in $D$ connecting the two cone points. After doing this for all such $D$, we get a set $C$. Isotop $g$ in $D$ so that $g(C) = C$. The components of $X - \mathrm{Int}N(C)$ are homeomorphic to the hyperbolic components of $X - \mathrm{Int}N(C')$. Therefore the modified $g$ is a reducible map on $X$ with $C$ the reducing set.   $\square$

**Theorem 1.5.** *Let $\widetilde{X}$ be a compact covering surface of a hyperbolic orbifold $X$. If $f : X \to X$ is a pseudo-Anosov map and $\widetilde{f} : \widetilde{X} \to \widetilde{X}$ is a lifting of $f$, then $\widetilde{f}$ is a pseudo-Anosov map in the sense of Thurston.*

*Proof.* Let $Y$ be the underlying surface of $X$. Then $f$ defines a generalized pseudo-Anosov map on $Y$, in the sense that the cone points might be singular points of the foliations with prong number $p = 1$. Since $\widetilde{X}$ is a surface, if $y$ is a point in $\widetilde{X}$ which covers a cone point of angle $2\pi/k$, then the foliations have prong number $kp \geq 2$ at $y$. It follows that the liftings of the foliations are foliations with no prong number 1 singularities in the interior of $\widetilde{X}$, and are preserved by $\widetilde{f}$. Therefore $\widetilde{f}$ is a pseudo-Anosov map in the sense of Thurston (see [FLP]).   $\square$

Suppose that $S$ is a hyperbolic surface. Let $\Lambda$ be a finite group acting on $S$. Then $X = S/\Lambda$ is a hyperbolic orbifold. If $f : S \to S$ is a $\Lambda$-equivariant map, then it induces a map $f' : X \to X$. Thus the theory of equivariant maps on surfaces follows from that on orbifolds. To be more precise, we make the following definitions.

**Definition 1.6.** (1) Let $\Lambda$ be a finite group acting on a hyperbolic surface $S$. A map $f : S \to S$ is a $\Lambda$-map if $f \circ \lambda = \lambda \circ f$ for all $\lambda \in \Lambda$.

(2) An isotopy $f_t$, $0 \leq t \leq 1$, is a $\Lambda$-isotopy if $f_t$ is a $\Lambda$-map for all $t$.

(3) A $\Lambda$-map $f : S \to S$ is a $\Lambda$-reducible map if there is an admissible set $C = C_1 \cup \ldots \cup C_n$ such that $f(C) = C$, and $\lambda(C) = C$ for all $\lambda \in \Lambda$.

(4) A $\Lambda$-map $f : S \to S$ is a $\Lambda$-pseudo-Anosov map if $f$ is a pseudo-Anosov map in the sense of Thurston, and the measured foliations are $\Lambda$-invariant.

Since any $\Lambda$-map $f$ on $S$ induces a map $f'$ on the orbifold, and an isotopy of $f'$ induces a $\Lambda$-isotopy of $f$, the following $\Lambda$-equivariant version of Thurston classification theorem is a direct consequence of Theorem 1.4.

**Corollary 1.7.** *If $f : S \to S$ is a $\Lambda$-map, then it is $\Lambda$-isotopic to a $\Lambda$-map $g$ which is either a periodic $\Lambda$-map, a $\Lambda$-reducible map, or a $\Lambda$-pseudo-Anosov map.*   $\square$

In [Be] Bers proved Thurston's classification theorem using quasi-conformal mapping theory. Corollary 1.7 is also true in this category. A conformal structure $\sigma$ on $S$ is $\Lambda$-invariant if all maps in $\Lambda$ are conformal maps with respect to $\sigma$. Similarly, a quadratic differential $\omega$ is $\Lambda$-invariant if

$$\lambda^*(\omega) = \omega \qquad \text{or} \qquad \lambda^*(\omega) = \overline{\omega}$$

depending on whether $\lambda$ is orientation preserving or orientation reversing. A $\Lambda$-map $f$ is a $\Lambda$-pseudo-Anosov map in the sense of Bers if there is a $\Lambda$-invariant conformal structure $\sigma$ on $S$, such that $f$ is an extremal map as defined in [Be], and its initial and terminal differential $\omega$ is $\Lambda$-invariant. Notice that such a differential defines a pair of $\Lambda$-invariant transverse measured foliations, preserved by the map $f$. Hence such a map is a $\Lambda$-pseudo-Anosov map as defined before. The proof of [Wu] can be modified to prove Corollary 1.7 in such stronger sense.

**Definition 1.8.** A reducing set $C$ of $f$ is called a *complete reducing set* if $f$ is isotopic to a map $g$ with $g(N(C)) = N(C)$, and $g^{n_i}$ restricted to $X_i$ is either a periodic map or a pseudo-Anosov map, where $X_i$ is an arbitrary component of $X - \text{Int}N(C)$, and $n_i$ is the least positive integer such that $g^{n_i}$ maps $X_i$ to itself.

**Remark 1.9.** Clearly a maximal reducing set of $f$ is a complete reducing set. A minimal complete reducing set is called a *canonical reducing set*. It can be shown that if $C$ is a canonical reducing set of $f$, then its lifting $\widetilde{C}$ is a canonical reducing set of a lifting of $f$ on a compact covering surface $\widetilde{X}$ of $X$. Since the canonical reducing set of a surface map is unique [Wu, Theorem 1], it follows that the canonical reducing set for $f$ is also unique, and it is empty if and only if $f$ is isotopic to a periodic or pseudo-Anosov map.

We now consider a fiber preserving map $f$ on an orientable Seifert fiber space $M$ which has hyperbolic orbifold $X(M)$. Let $p : M \to X(M)$ be the projection map. Let $\widehat{f} : X(M) \to X(M)$ be the orbifold map induced by $f$. Recall that a torus $T$ in $M$ is a vertical torus if it is a union of fibers in $M$.

**Lemma 1.10.** *There is a collection of vertical tori $\mathcal{T}$ in $M$, and a map $g$ which is isotopic to $f$ via a fiber preserving isotopy, such that*

(1) $g(N(\mathcal{T})) = N(\mathcal{T})$;

(2) *Each component of $M - \text{Int}N(\mathcal{T})$ either is a twisted $I$-bundle over Klein bottle, or has hyperbolic orbifold; and*

*(3) If $M_i$ is a component of $M - IntN(\mathcal{T})$ with hyperbolic orbifold, and $g(M_i) = M_i$, then the orbifold map $\widehat{g} : X(M_i) \to X(M_i)$ is either periodic or pseudo-Anosov.*

*Proof.* Recall that since $X(M)$ is the orbifold of an orientable Seifert fiber space, it has no reflectors [Sc]. Choose a complete reducing set $C$ on $X(M)$ for the map $\widehat{f}$, and modify it as follows: If $\alpha$ is an arc in $C$ connecting two cone points of order 2, replace it by the curve $\alpha' = \partial N(\alpha)$. We can modify $f$ via a fiber preserving isotopy, to get a map $g$ such that $\widehat{g}(N(C)) = N(C)$, and the restriction of $\widehat{g}$ to each hyperbolic component of $X(M) - \text{Int}N(C)$ is either periodic or pseudo-Anosov. It is now clear that the map $g$ and the collection of tori $\mathcal{T} = p^{-1}(C)$ satisfy the conclusion of the lemma.   $\square$

## §2. Classification of maps on Euclidean orbifolds

A euclidean orbifold is an orbifold $X$ such that its orbifold Euler number is zero. A curve $C$ in $X - S(X)$ is an *essential curve* of $X$ if no component of $X - \text{Int}N(C)$ is a disk with at most one cone point, or an annulus containing no singular points of $X$ and having one boundary on $\partial X$. Recall that a map $f$ on a torus is an *Anosov map* if there are a pair of transverse measured foliations $(\mathcal{F}_1, u_1)$ and $(\mathcal{F}_2, u_2)$ such that $f$ preserves the foliations, and the measures $u_i$ are changed by $f(u_1) = \lambda u_1$, and $f(u_2) = (1/\lambda)u_2$, where $\lambda > 1$. We do not define Anosov maps on other euclidean *surfaces*: as we will see, all maps on them are isotopic to periodic ones.

**Definition 2.1.** Let $f : X \to X$ be an orbifold map on a euclidean orbifold $X$.

(1) $f$ is *periodic* if $f^n = id$ for some $n$.

(2) $f$ is *reducible* if $f(C) = C$ for some essential curve $C$ of $X$. We call $C$ a *reducing curve* of $f$.

(3) If $S(X) \neq \emptyset$, then $f : X \to X$ is an *Anosov* map if $f : X' \to X'$ is a pseudo-Anosov map in the sense of Thurston, where $X' = X - S(X)$.

**Lemma 2.2.** *Any map $f : K \to K$ on a Klein bottle is isotopic to a periodic map.*

*Proof.* The fundamental group and the first homology group of $K$ have the following presentation

$$\pi_1(K) = \left\langle a, b \,|\, aba^{-1} = b^{-1} \right\rangle,$$

$$H_1(K) = \mathbb{Z} \oplus \mathbb{Z}_2 = \left\langle a, b \,|\, ab = ba, b^2 = 1 \right\rangle,$$

where $a$ and $b$ are represented by simple closed curves on $K$. Any element of $\pi_1(K)$ can be written uniquely as $b^m a^n$ in $\pi_1(K)$, which represents $mb + na$ in $H_1(K)$. So if $C$ is an oriented simple closed curve on $K$ representing $b$ in $H_1(K)$, then $n = 0$ and $m \cong 1 \mod 2$. Since $C$ is simple, we must have $m = \pm 1$. Hence $C$ is isotopic to $b^{\pm 1}$ as oriented curves. Since $f(C)$ also represents $b$ in $H_1(K)$, we must have $f_*(C) = C^{\pm 1}$ in $\pi_1(K)$. Notice that $aba^{-1} = b^{-1}$, so $aba^{-1}$ is freely homotopic to $b$, and hence $f(C)$ is isotopic to $C$ as oriented curves. By an isotopy we may assume that $f$ is an identity map on $C$, so we may further assume that $f^2$ is identity in $N(C)$. After cutting along $C$, we get an annulus $A$ such that $f^2|_{\partial A} = id$. Thus $f^2$ must be a twist on $A$. Notice that $C$ induces opposite orientations on the two components of $\partial A$, hence a $4\pi$ twist on $A$ is isotopic to the identity map on $K$. Thus $f^4$ is always isotopic to the identity map. By the Nielsen realization theorem for periodic maps on surfaces, $f$ is isotopic to a map $g$ such that $g^4 = id$.   □

If $F$ is a surface, we use $F(p_1, \ldots, p_k)$ to denote an orbifold with underlying surface $F$, and with cone points of order $p_1, \ldots, p_k$. If we use $S, P, T, K, D, A, U$ to denote sphere, projective plane, torus, Klein bottle, disk, annulus and Möbious band, respectively, then a euclidean orbifold without reflector is one of the following:

(1) $T$, $K$, $A$, or $U$;

(2) $D(2, 2)$;

(3) $P(2, 2)$;

(4) $S(2, 3, 6)$, $S(2, 4, 4)$, $S(3, 3, 3)$;

(5) $S(2, 2, 2, 2)$.

**Lemma 2.3.** *If $f : X \to X$ is a map on $X = P(2, 2)$, then $f$ is isotopic to a periodic map.*

*Proof.* Write $P(2, 2) = U(2, 2) \cup D$. By an isotopy we may assume that $f|_D = id$. Let $C$ be the central curve of $U$ containing the two cone points. Then it is clear that $f|_{U(2,2)}$ can be rel $\partial U \cup S(X)$ isotoped so that $f(C) = C$. By a further orbifold isotopy we may assume that $f$ is periodic in $N(C)$. Since $P(2, 2) - \text{Int} N(C)$ is a disk $D'$, $f$ can be rel $\partial D'$ isotoped to be periodic on $D'$ as well.   □

**Proposition 2.4.** *If a euclidean orbifold $X$ is not $T$ or $S(2, 2, 2, 2)$, then any map $f$ on $X$ is isotopic to a periodic one.*

*Proof.* First assume that $X$ has no reflectors. The result is obvious when $X = A$ or $U$, and has been proved for $K$ and $P(2, 2)$ in Lemmas 2.2 and 2.3 respectively. In all other cases $X - S(X)$ is a pair of pants, so any map on $X - S(X)$ (hence on $X$) is periodic.

Now assume $X$ has no reflector circles but has some reflector arcs. By an isotopy we may assume that $f$ is already a periodic map on $N(C)$, where $C$ is the union of $\partial X$ and the reflectors. Since $X$ has no reflector circles but some reflector arcs, the Euler characteristic of $N(C)$ is negative, no matter whether it has reflector corners or not. Thus $X - \operatorname{Int} N(C)$ is an orbifold with nonempty boundary, which has positive Euler number but no reflectors, hence it must be a disk with at most one cone point. It follows that $f$ is rel $N(C)$ isotopic to a periodic map.

If $X$ has some reflector circles, one can ignore these reflector circles to get another orbifold $X'$. Clearly $f : X \to X$ is isotopic to a periodic map if and only if $f : X' \to X'$ is, so the result follows from the above two cases.  $\square$

**Theorem 2.5.** *(a) A map $f$ on a euclidean orbifold $X$ is isotopic to either a periodic map, a reducible map, or an Anosov map.*

*(b) Let $\widetilde{X}$ be a surface that covers $X$, and let $\widetilde{f} : \widetilde{X} \to \widetilde{X}$ be a lifting of $f$. If $f$ is periodic (resp. reducible, Anosov) then $\widetilde{f}$ is also periodic (resp. reducible, Anosov).*

*Proof.* (a) If $X = T$, this is a well known classical result. If $X = S(2,2,2,2)$, then $X - S(X)$ is a hyperbolic surface, so the result follows immediately from the definition above and Thurston's classification theorem. In all other cases by Proposition 2.4 the map $f$ is isotopic to a periodic map.

(b) This is quite obvious if $f$ is periodic or reducible, or if $X = T$. Suppose $X = S(2,2,2,2)$, $\widetilde{X} = T$, and $f$ is Anosov. By the proof of Theorem 1.5, the measured foliations on $X$ lifts to measured foliations on $\widetilde{X}$ which has no prong number 1 singularities. Recall that if a foliation on a surface $F$ has singular points of prong number $p_1, \ldots, p_k$, then the Euler characteristic of $F$ is $\sum (2 - p_i)/2$. Since the Euler characteristic of $\widetilde{X}$ is zero, the lifted foliations on $\widetilde{X}$ have no singularities. Thus by definition $\widetilde{f}$ is Anosov. The proof in the case that $\widetilde{X} = X = S(2,2,2,2)$ is similar.  $\square$

**Remark 2.6.** Not all results of hyperbolic orbifolds generalize to euclidean orbifolds. For example, if $X$ is a hyperbolic orbifold, then two maps on $X$ are isotopic if and only if their liftings are isotopic. However, there is a map on $S(2,2,2,2)$ which permutes the singular points (so it is not isotopic to the identity map), and yet its lifting to the torus is isotopic to the identity map.

§3. INVARIANT SEIFERT FIBRATIONS FOR MAPS

In this section we will study the problem of which manifolds $M$ have a Seifert fibration such that all orientation preserving maps on $M$ are isotopic to a fiber preserving map with respect to this fibration. Theorems 3.1 and 3.2 are what we can find in the literature. The results in this section will deal with the remaining cases. It turns out that all but four Seifert fiber spaces have this property. Moreover, if the manifold is not a lens space or $S^3$, then the fibration is preserved by orientation reversing maps as well. See Theorem 3.11 below.

**Theorem 3.1.** ([Sc, Theorem 3.9]) *Let $M, N$ be compact Seifert fiber spaces (with fixed Seifert fibrations), and let $f : M \to N$ be a homeomorphism. Then $f$ is homotopic to a fiber preserving homeomorphism (and hence an isomorphism of Seifert bundles) unless one of the following occurs.*

*(a) $M$ is covered by $S^3$ or $S^2 \times \mathbb{R}$,*

*(b) $M$ is covered by $T^3$,*

*(c) $M$ is $S^1 \times D^2$ or a I-bundle over the torus or Klein bottle.* □

It is now known that two homeomorphisms of Seifert fiber spaces are homotopic if and only if they are isotopic, see [Wa] for Haken manifolds, [Sc2] for irreducible Seifert manifolds with infinite fundamental groups, and [BO, RB, HR, R, B, La] for various cases of Seifert manifolds covered by $S^3$ and $S^2 \times S^1$. Thus the word "homotopic" in the theorem can be replaced by "isotopic". The theorem says that if $M$ is not one of the listed, then Seifert fibrations on $M$ are unique up to isotopy. For $S^1 \times D^2$ and $T \times I$, there are obviously maps which is not isotopic to any fiber preserving maps with respect to any Seifert fibration. If $M$ is a twisted $I$-bundle over the Klein bottle, then $M$ has two Seifert fibrations up to isotopy [WW, Lemma 1.1]. Note that those two Seifert fibrations have different orbifolds. Let $p : M \to X(M)$ be one of the Seifert fibrations. Then $p \circ f : M \to X(M)$ is another Seifert fibration with the same orbifold. Since these two Seifert fibrations have the same orbifold, [WW, Lemma 1.1] implies that they are isotopic, that is, there is a map $\varphi$ of $M$, isotopic to the identity, so that $p \circ f \circ \varphi = p$. It follows that $f$ is isotopic to a fiber preserving map with respect to each of the two Seifert fibrations of $M$. It remains to consider the manifolds covered by $S^3$, $S^2 \times \mathbb{R}$, or $T^3$.

**Theorem 3.2.** ([BO, Proposition 3.1]) *Suppose $M$ is a Seifert manifold with orbifold $S(2, 3, p)$, $p \geq 5$ or $S(3, 3, q)$, $q \geq 2$. Then all homeomorphisms $f$ of $M$ are isotopic to a fiber preserving map of $M$.* □

By [Ja, VI.16], Seifert fibration structures on the manifolds in the theorem are unique up to isomorphism. The above theorem says that it is unique up to isotopy.

In this section we will discuss the remaining cases. The following lemma is useful in detecting maps which are isotopic to fiber preserving maps. Note that the lemma is not true if $X(M)$ is a $S(p,q)$ or $D(p)$, since there are many non isotopic fibrations on lens spaces and the solid torus, all having one fiber in common.

**Lemma 3.3.** *Let $p : M \to X(M)$ be a Seifert fibration of a close orientable 3-manifold $M$, such that $X(M)$ is not a sphere with at most two cone points or a disk with at most one cone point. Let $f : M \to M$ be an homeomorphism. If there is a fiber $\alpha$ of $M$ such that $f(\alpha)$ is isotopic to another fiber, then $f$ is isotopic to a fiber preserving map with respect to the fibration $p$.*

*Proof.* By an isotopy we may assume that $f(\alpha)$ is a fiber of $M$, and $f$ maps a fibered neighborhood $N(\tau)$ to a fibered neighborhood $N(f(\alpha))$. Let $M' = M - \mathrm{Int}N(\alpha)$, and let $N' = M - \mathrm{Int}N(f(\alpha))$, with Seifert fibrations inherited from that of $M$. These are manifolds with boundary, so by Theorem 3.1 the map $f|_{M'}$ is isotopic to a fiber preserving map unless $M'$ is one of the manifolds listed in Theorem 3.1(c). If $M'$ is a solid torus or a $T \times I$, then $X(M)$ is either a sphere with at most two cone points or a disk with at most one cone point, which have been excluded in the lemma. The remaining case is that $M'$ and $N'$ are twisted $I$-bundles over the Klein bottle. Recall that up to isotopy there are only two Seifert fibrations on a twisted $I$-bundle over Klein bottle, with orbifold $D(2,2)$ and the Möbius band respectively [WW, Lemma 1.1]. Now the orbifolds $X(M')$ and $X(N')$ are both obtained from $X(M)$ by removing a disk with possibly one cone, so they are both orientable or both non-orientable. It follows that $f|_{M'}$ is isotopic to a fiber preserving map. Extending the isotopy over $M$, we get a map $g : M \to M$, which is fiber preserving on $M'$. Now $g|_{N(\alpha)}$ maps the fibered solid torus $N(\alpha)$ to a fibered solid torus $N(f(\alpha))$, sending a fiber on $\partial N(\alpha)$ to a fiber on $\partial N(f(\alpha))$. By [Ja, Lemma VI.19], $g|_{N(\alpha)}$ is rel $\partial N(\alpha)$ isotopic to a fiber preserving map.   $\square$

Denote by $\mathcal{H}(M)$ the homeotopy group of $M$, i.e., the group of all homeomorphisms modulo the subgroup of all those isotopic to identity. Denote by $\mathcal{H}^+(M)$ the subgroup represented by orientation preserving maps. We do not consider $S^3$ or $S^2 \times S^1$ as lens spaces.

**Lemma 3.4.** *Any lens space $M$ has a Seifert fibration $p : M \to X(M)$, such that any orientation preserving maps on $M$ is isotopic to a fiber preserving map with respect to this fibration.*

*Proof.* According to Hodgson and Rubinstein [HR, Theorem 5.6], $\mathcal{H}(M) = \mathbb{Z}_2$,

$\mathbb{Z}_2 + \mathbb{Z}_2$, or $\mathbb{Z}_4$. Moreover, by [HR, Corollary 5.7], if $\mathcal{H}(M) = \mathbb{Z}_4$ then its two generators are represented by orientation reversing maps. All other elements of $\mathcal{H}(M)$ are represented by involutions which preserve a Heegaard splitting torus $T$ in $M$. Hence up to isotopy we may assume that $f(T) = T$, and $f$ is an involution. The torus $T$ cuts $M$ into two solid tori $V_1$ and $V_2$. If $f$ maps each $V_i$ to itself, then $f|_{V_i}$ is an orientation preserving involution, which preserves any Seifert fibration of $V_i$ up to isotopy. So we may assume that $f$ maps $V_1$ to $V_2$. Then $f|_T$ is an involution. If $f|_T$ has a fixed point, then its fixed point set is a union of 2 circles. If $f|_T$ has no fixed points, then $T/f$ is a Klein bottle, so any orientation reversing simple closed curve on it lifts to a curve $C$ on $T$ preserved by $f$. In any case, there is a simple closed curve $C$ on $T$ such that $f(C) = C$. The curve $C$ can not bound a disk on $V_i$, otherwise the properties $f(C) = C$ and $f(V_1) = V_2$ would imply that it bounds disks on both sides, so $M = S^2 \times S^1$, a contradiction. Construct an $S^1$ bundle structure on $T$ which is invariant under the involution $f$, with $C$ as a fiber, extend it to a Seifert fibration over one of the solid torus, then map it to the other side using $f$ to get Seifert fibration on the other side. Clearly, $f$ preserves this Seifert fibration of $M$.  $\square$

**Lemma 3.5.** *If $p : M \to X(M)$ is a Seifert fibration with $X(M) = S(2, 2, k)$ for some $k > 1$, then any map $f : M \to M$ is isotopic to a fiber preserving map with respect to this fibration.*

*Proof.* Let $\alpha$ be an arc on $X(M)$ connecting two cone points of order 2. Then $p^{-1}(\alpha)$ is a Klein bottle $K_1$ in $M$. Note that if $k = 2$ then there are three cone points of order 2, so there are three such arcs, and hence three Klein bottles $K_1, K_2, K_3$. It was shown in the proof of [R, Theorem 6] that $f$ can be isotoped so that if $k \neq 2$ then $f$ maps $K_1$ to itself, and if $k = 2$ then $f$ maps $K_1$ to one of the $K_i$. Note that each $K_i$ is a union of fibers of $M$, so it gives a fibration of $K_i$ with orbifold an arc with two cone points. A regular fiber of $K_i$ is a separating essential simple closed curve, which is unique up to isotopy, hence $f$ can be isotoped, by a global isotopy, so that $f$ preserves the fibration on $K_i$. Since the regular fibers of $K_i$ are regular fibers of $M$, it follows from Lemma 3.3 that $f$ is isotopic to a fiber preserving map.  $\square$

**Lemma 3.6.** *Suppose $p : M \to X(M)$ is a Seifert fibration of a spherical manifold $M$ with orbifold $X(M) = S(2, 3, 4)$. Then any map $f : M \to M$ is isotopic to a fiber preserving map with respect to this fibration.*

*Proof.* Identify $X(M)$ with the unit sphere in $\mathbb{R}^3$ so that the three cone points lie

on the $xy$-plane. Let

$$\widehat{g} : X(M) \to X(M)$$

be the reflection along the $xy$-plane. There is a fiber preserving involution

$$g : M \to M$$

with $\widehat{g}$ as its orbifold map, constructed as follows: Removing a small neighborhood $N(c_i)$ for each cone point $c_i$, we get a three-punctured sphere $Q \subset X(M)$, and $p^{-1}(Q) = Q \times S^1$. Define

$$g : Q \times S^1 \to Q \times S^1$$

by $g(x, y) = (\widehat{g}(x), \bar{y})$, where $\bar{y}$ is the complex conjugate of $y$. The restriction of $g$ on each boundary torus of $Q \times S^1$ is an orientation preserving involution, which preserves the slope of any curve. Hence it extends to a fiber preserving involution on each fibered solid torus $p^{-1}(N(c_i))$, resulting the required map $g$.

Let $P$ be a sphere with 3 open disks removed. The manifold $P \times S^1$ has a product fibration. Its fiber and the boundaries of $P$ set up a longitude-meridian pair $(l_i, m_i)$ for each boundary component $T_i$ of $P \times S^1$. Denote by $M(p_1, q_1; p_2, q_2; p_3, q_3)$ the manifold obtained by attaching a solid torus $V_i$ to each $T_i$, so that a meridian of $V_i$ is glued to the curve $p_i l_i + q_i m_i$ on $T_i$. The manifold $M$ can be written as $M(2, 1; 3, b_2; 4, b_3)$, where $b_2 = \pm 1$. By reversing the orientation of $M$ if necessary, we may assume that $b_2 = 1$, so $M = M(2, 1; 3, 1; 4, b_3)$. By [Or, 6.2 Theorem 2], the first homology group of $M$ is

$$H_1(M) = \mathbb{Z}_{|12 + 6b_3 + 8b_2|} = \mathbb{Z}_{|20 + 6b_3|},$$

hence $H_1(M) = \mathbb{Z}_2$ if and only if $b_3 = -3$. The first homology group of $M$ has a presentation

$$H_1(M) = \mathbb{Z}[x, y, z, t] / \langle 2x + t = 0, 3y + t = 0, 4z + b_3 t = 0 \rangle.$$

The map $g$ defined above sends each generator $u$ to $-u$. Hence if $H_1(M) = \mathbb{Z}_p$ and $p \neq 2$, then $g_* : H_1(M) \to H_1(M)$ is nontrivial. It follows that $g_*$ is nontrivial unless $b_3 = -3$.

Let $\mathcal{H}(M)$ be the group of all homeomorphisms module those isotopic to the identity. According to Rubinstein and Birman [RB, Main Theorem], $\mathcal{H}(M) = \mathbb{Z}_2$ or a trivial group, and it is trivial if and only if $b_3 = -3$. The above argument shows that when $b_3 \neq -3$, the map $g$ is a nontrivial element in $\mathcal{H}(M)$, hence any map on $M$ is isotopic to either $g$ or the identity map. When $b_3 = -3$, $\mathcal{H}(M)$ is trivial, so all maps are isotopic to the identity map. Since both $g$ and the identity map are fiber preserving, the result follows. $\square$

**Lemma 3.7.** *Let $\rho : M \to S^1$ be a torus bundle structure on $M$. Let $f : M \to M$ be a homeomorphism. Then either $M$ is an $S^1$ bundle over $T^2$, or $f$ is isotopic to a torus fiber preserving map with respect to the fibration $\rho$.*

*Proof.* Let $T$ be a torus fiber in $M$, and let $T' = f(T)$. Isotop $f$ so that $|T \cap T'|$ is minimal. After cutting along $T$, we get a manifold $N = T \times I$, so if $|T \cap T'| = 0$, we can isotop $f$ so that $f(T) = T$, and a further isotopy will deform $f$ to a torus fiber preserving map, and we are done. So we assume that $|T \cap T'| \neq 0$.

By an innermost circle argument one can eliminate all trivial circles in $T \cap T'$. Thus each component of $T' - T$ is an annulus. Also, a boundary compression of an annulus of $T' - T$ would isotop $T'$ to a torus having less intersection with $T$. Thus the minimality of $|T \cap T'|$ implies that these annuli are $\partial$-incompressible in the manifold $N$ obtained by cutting $M$ along $T$, so they are essential annuli $P$ in $N$. There is an $S^1$ bundle structure of $N$ over an annulus $A$, so that $P$ is a union of fibers. $M$ can be recovered by gluing the two sides of $\partial N$ together. Since the boundary of $P$ matches each other to form the torus $T'$, the $S^1$ bundle structure extends over $M$, whose orbifold $X(M)$ is obtained by gluing the two boundary components of $A$ together, hence is a torus. $\square$

If $M$ is a union of two twisted $I$-bundles over Klein bottle glued together along their boundary, then there is a map

$$\rho : M \to J = [-1, 1],$$

such that $\rho^{-1}(x)$ is a Klein bottle if $x = \pm 1$, and is a torus otherwise. The map $\rho$ is called a torus bundle structure of $M$ over $J$, while $J$ is considered as a closed 1-orbifold with two singular points at $\pm 1$. The proof of the following lemma is similar to that of Lemma 3.7. We omit the details.

**Lemma 3.8.** *Let $\rho : M \to J$ be a torus bundle structure of $M$ over $J$. Let $f : M \to M$ be a homeomorphism. Then either $M$ is a Seifert fiber space over a Euclidean orbifold, or $f$ is isotopic to a torus fiber preserving map with respect to the fibration $\rho$.* $\square$

**Lemma 3.9.** *Suppose $p : M \to X(M)$ is a Seifert fibration of a Euclidean manifold $M$ with orbifold $X(M) = S(2, 4, 4)$ or $S(2, 2, 2, 2)$. Then any map $f : M \to M$ is isotopic to a fiber preserving map with respect to this fibration.*

*Proof.* First assume $X(M) = S(2, 4, 4)$. Put $T = S^1 \times S^1$, with each point represented by $(x, y)$, where $x, y$ are mod 1 real numbers. Then the manifold $M$ can be expressed as $(T \times I)/\varphi$, with gluing map

$$\varphi : T \times 0 \to T \times 1$$

given by $\varphi(x, y, 0) = (-y, x, 1)$. Consider the universal covering $\widetilde{W}$ of $W = T \times I$, which is identified with $\mathbb{R}^2 \times I$. The lifting of the two order 4 singular fibers of $M$ in $\widetilde{W}$ are the vertical arcs at $(x, y)$, with both $x, y$ integers, or both half integers.

CLAIM 1. *There is a fiber preserving map $f'$ on $M$ which reverses the fiber orientation.*

On $W$ the map is given by $f'(x, y, z) = (y, x, 1 - z)$. Since

$$f' \circ \varphi(x, y, 0) = f'(-y, x, 1) = (x, -y, 0) = \varphi^{-1} \circ f'(x, y, 0),$$

the map $f'$ induces a map on $M$, which clearly is fiber preserving and fiber orientation reversing. The claim follows.

The product structure $W = T \times I$ gives a torus bundle structure of $M$ over $S^1$. By Lemma 3.7 we may assume that $f$ preserves the torus fibration. Note that the map $f'$ constructed in the claim is also torus fiber preserving, and it reverses the orientation of $S^1$, the torus bundle orbifold. If $f$ also reverses the orientation of $S^1$, then $f' \circ f$ is orientation preserving on $S^1$, and $f$ is isotopic to a Seifert fiber preserving map if and only if $f' \circ f$ is. Therefore, by considering $f' \circ f$ instead if necessary, we may assume without loss of generality that $f$ preserves the orientation of $S^1$. After an isotopy, we may further assume that $f$ induces an identity map on $S^1$, so it maps each horizontal torus $T \times t$ to itself. We say that $f : M \to M$ is a level preserving map.

Let $\alpha : I \to M$ be given by $\alpha(t) = (0, 0, t)$. It is a path representing a singular fiber. By an isotopy we may assume that $f$ fixes the base point of $\alpha$. Let $\alpha'$ be the path $f \circ \alpha$.

CLAIM 2. *The path $\alpha'$ is isotopic to a singular fiber in $M$.*

Let $\widetilde{\alpha}$ and $\widetilde{\alpha}'$ be the lifting of $\alpha$ and $\alpha'$ on $\widetilde{W}$ with initial point $(0, 0, 0)$. Let $\widetilde{\beta}'$ be the straight line in $\widetilde{W}$ connecting the ends of $\widetilde{\alpha}'$. Since $f$ is level preserving, $\widetilde{\alpha}'$ intersect each plane $\mathbb{R}^2 \times t$ at a single point, so there is an isotopy $\psi_t$ from $\widetilde{\alpha}'$ to $\widetilde{\beta}'$ rel $\partial$, moving each point of $\widetilde{\alpha}'$ along the straight line on a level plane to the corresponding point of $\widetilde{\beta}'$. For each $t$, the image of the path $\psi_t$ under the projection $q : \widetilde{W} \to M$ is a closed curve, which has no self intersection because $\psi_t$ intersects each level plane at a single point. Hence $\alpha'$ is isotopic to the curve $\beta' = q \circ \widetilde{\beta}'$. Without loss of generality, we may assume that $\widetilde{\alpha}'$ is a straight line, connecting $(0, 0, 0)$ to $(x, y, 1)$. Note that $x, y$ must be integers.

Put

$$(\bar{x}, \bar{y}) = (\frac{x - y}{2}, \frac{x + y}{2}).$$

For any $t$, let

$$h_t : I \to \mathbb{R}^2 \times I$$

be the straight line from $t(\bar{x}, \bar{y}, 0)$ to the point $(1 - t)(x, y, 1) + t(\bar{x}, \bar{y}, 1)$. Each $h_t$ projects to a simple arc in $T \times I$. By definition, we have

$$\begin{aligned}
\varphi(h_t(0)) &= \varphi(t(\bar{x}, \bar{y}, 0)) = t(-\bar{y}, \bar{x}, 1) \\
&= (\frac{t(-x - y)}{2}, \frac{t(x - y)}{2}, 1) \\
&= (1 - t)(x, y, 1) + t(\bar{x}, \bar{y}, 1) \\
&= h_t(1) \qquad \mathrm{mod}\ 1,
\end{aligned}$$

so each $h_t$ gives a closed curve in $M$. Moreover, since $h_t$ intersects each level torus at a single point, it has no self intersection. Thus $h_t$ is an isotopy of $\alpha'$. When $t = 1$, the loop $h_1$ is a vertical line $h_1(s) = (\bar{x}, \bar{y}, s)$, which is a singular fiber of $M$ because $(\bar{x}, \bar{y})$ are both integers or both half integers. This proves the claim. Figure 3.1 shows the isotopy in $\mathcal{R}^2 \times I$, when $(x, y) = (1, 0)$. The vertical line at $(0, 0, 0)$ is the curve $\widetilde{\alpha}$, the line from $(0, 0, 0)$ to $(1, 0, 1)$ is $\widetilde{\alpha}$, the vertical line at $(\bar{x}, \bar{y}, 0) = (\frac{1}{2}, \frac{1}{2}, 0)$ is the image of $h_1$, which is the lift of a singular fiber in $M$, and the lines between them represent the isotopy $h_t$.

The isotopy extends to an isotopy of $M$, modifying $f$ to a map $g$, which maps the singular fiber $\alpha$ to another singular fiber. By Lemma 3.3 such a map is isotopic to a fiber preserving map.

The proof for the case $X(M) = S(2, 2, 2, 2)$ is similar. The gluing map $\varphi : T \times 0 \to T \times 1$ is now given by $\varphi(x, y, 0) = (-x, -y, 1)$. Write $(\bar{x}, \bar{y}) = (x/2, y/2)$, and define $h_t$ using the same formula as above. The remaining part of the proof are the same, only notice that in this case a vertical line of $\widetilde{W}$ at $(x, y)$ represents a singular fiber of $M$ if and only if each of $x$ and $y$ is an integer or a half integer. $\quad \square$

Figure 3.1

We now consider the manifold $M = S^2 \times S^1$. Put $S^2 = D_1^2 \cup D_2^2$, and $V_i = D_i^2 \times S^1$, so $M$ is the union of $V_1$ and $V_2$ along a torus $T$. Let $m, l$ be a meridian-longitude pair on $T$, where $m$ bounds a disk in both $V_i$. The torus $T$ is a Heegaard splitting torus of $M$, which is unique up to isotopy [Wa2]. Thus any map $f : M \to M$ can be isotoped so that $f(T) = T$. The map $f|_T$ corresponds to a matrix

$$\begin{pmatrix} p & q \\ r & s \end{pmatrix}$$

so that $f_*([l]) = p[l] + q[m]$ and $f_*([m]) = r[l] + s[m]$ in $H_1(T)$. Since $f|_T$ extends over the solid tori, we have $r = 0$ and $s = \pm 1$, so $p = \pm 1$. Note that any map on $T$ extends to a unique orientation preserving map on $M$ up to isotopy, which maps each $V_i$ to itself if and only if the determinant of the matrix equals 1.

Denote by $\alpha, \beta, \gamma$ the orientation preserving maps on $M$ corresponding to the following matrices, respectively.

$$\begin{pmatrix} -1 & 0 \\ 0 & -1 \end{pmatrix}, \begin{pmatrix} 1 & 0 \\ 0 & -1 \end{pmatrix}, \begin{pmatrix} 1 & 1 \\ 0 & -1 \end{pmatrix}$$

These matrices generate all the matrices satisfying the conditions that $r = 0$ and $s = \pm 1$, hence the maps $\alpha, \beta, \gamma$ generate the group $\mathcal{H}^+(M)$ of all isotopy classes of orientation preserving maps on $M$, which is an index 2 subgroup of the homeotopy group $\mathcal{H}(M)$.

There is an isotopy of $S^2$ which changes $m$ to $-m$, the product of which with the identity map on $S^1$ then isotop $\beta$ to the identity map, so $\beta$ represents the identity in $\mathcal{H}^+(M)$. Also, $\alpha$ and $\gamma$ both have order 2, and they commute, so they generate a group $\mathbb{Z}_2 + \mathbb{Z}_2$. It is known that this is indeed the group $\mathcal{H}^+(M)$.

Each non meridional curve on $T$ extends to a unique Seifert fibration of $M$. Consider the Seifert fibration corresponding to the curve $C$ which represents $[m] + 2[l]$ in $H_1(T)$. Clearly, $\alpha$ maps the isotopy class of $C$ to itself, with orientation reversed. Since

$$\gamma_*[C] = \gamma_*[m] + 2\gamma_*[l] = -[m] + 2([l] + [m]) = [C],$$

we see that $\gamma$ also maps the isotopy class of $C$ to itself. Thus they both preserves the Seifert fibration corresponding to $C$. There are infinitely many isotopy classes of Seifert fibrations on $S^2 \times S^1$, all of which are preserved by the map $\alpha$, but it is interesting that the fibration above is the only one preserved by $\gamma$. Let $\delta$ be the orientation reversing involution of $M$ such that $\delta|_T = id$. Then $\delta$ also preserves all curves on $T$, hence all Seifert fibrations of $M$. Since $\alpha, \gamma, \delta$ generate the homeotopy group of $M$, it follows that the fibration corresponding to $C$ is preserved by all maps on $M$. We have proved the following lemma for $S^2 \times S^1$. Recall from [Sc] that the two manifolds listed are the only ones which admit the $S^2 \times \mathbb{R}$ geometry.

**Lemma 3.10.** *Suppose $M$ is either $S^2 \times S^1$ or $RP^3 \# RP^3$. Then there is a Seifert fibration on $M$ such that any orientation preserving map on $M$ is isotopic to a fiber preserving map with respect to this fibration.*

*Proof.* We have seen that the lemma is true for $S^2 \times S^1$, so we assume that $M = RP^3 \# RP^3$. Let

$$\varphi = \varphi_1 \times \varphi_2 : S^2 \times S^1 \to S^2 \times S^1,$$

where $\varphi_1$ is an antipodal map on $S^2$, and $\varphi_2 : S^1 \to S^1$ is given by the involution $\varphi_2(x, y) = (-x, y)$, with $S^1$ considered the unit circle in $R^2$. Then the manifold $M$ is a quotient

$$M = (S^2 \times S^1)/\varphi.$$

Let $\eta = \varphi_1 \times id$, and $\psi = \varphi_1 \times \psi_2$, where $\psi_2(x, y) = (x, -y)$. Clearly, the maps $\eta$ and $\psi$ commute with $\varphi$, so they induce maps $\widehat{\eta}, \widehat{\psi} : M \to M$, which reverse the orientation of a reducing sphere in $M$. Moreover, the product $S^1$ bundle structure of $S^2 \times S^1$ is preserved by all the three maps, so it induces a Seifert fibration on $M$, which is preserved by $\widehat{\eta}$ and $\widehat{\psi}$. The map $\widehat{\eta}$ is orientation reversing. Thus to prove the lemma, it remains to show that any orientation preserving map on $M$ is isotopic to either $\widehat{\psi}$ or the identity.

Let $S$ be an essential 2-sphere in $M$, which is unique up to isotopy because $M$ has only two irreducible prime factors, so any map $f$ on $M$ can be isotoped so that $f(S) = S$. Since $\widehat{\psi}$ reverses the orientation of $S$, exactly half of the isotopy classes of orientation preserving maps on $M$ preserve the orientation of $S$. Therefore, we need only show that if $f$ preserves the orientations of both $M$ and $S$, then $f$ is isotopic to the identity map. Since any orientation preserving map on $S$ is isotopic to the identity map, we may assume that $f|_S = id$. The sphere $S$ cuts $M$ into two pieces $X_1, X_2$, each being a punctured $RP^3$, which is a twisted $I$-bundle over the projective plane $P = RP^2$. Since $f$ is orientation preserving, it maps each $X_i$ to itself. The result now follows from the following claim.

CLAIM. *If $X$ is a twisted $I$-bundle over the projective plane $P$, and $f : X \to X$ is an orientation preserving map such that $f|_{\partial X} = id$, then $f$ is rel $\partial X$ isotopic to an identity map.*

Clearly, the projective plane $P$ in $X$ is unique up to isotopy. We can first deform $f$ so that $f$ sends $P$ to itself, then further isotop $f$ so that the restriction of $f$ on $P$ is the identity.

Put $Y = S^2 \times [0,1]$. View $X$ as the quotient space of $(S^2 \times [0,1])/\eta$, where

$$\eta : S^2 \times 1 \to S^2 \times 1$$

identifies each point on the sphere $S^2 \times 1$ to its antipodal point. Let $\widetilde{f} : Y \to Y$ be the lift of $f$. Now the restriction of $\widetilde{f}$ on $S^2 \times 0$ is the identity, and the restriction of $\widetilde{f}$ on $S^2 \times 1$ is the lifting of the identity map on $P$, hence is either the identity or the antipodal map; but since it is the restriction of the orientation preserving map $\widetilde{f}$ on a boundary component of $Y$, it can not be the antipodal map. Thus the restriction of $\widetilde{f}$ on $\partial Y$ is the identity.

Consider $S^2$ as the standard unit sphere in $\mathbb{R}^3$, and let $x$ be the north pole. Then $C = x \times [0,1]$ is an arc in $Y$. By the Light Bulb Theorem [Ro, Page 257], we may isotopy $\widetilde{f}$ via an isotopy rel $\partial Y$, so that $\widetilde{f}|_C = id$. Let $D_1, D_2$ be the upper and lower half disks of $S^2$. Isotop $\widetilde{f}$ rel $\partial Y$ so that it maps each $D_1 \times t$ isometrically to itself, then further isotop $\widetilde{f}$ on $D_2 \times [0,1]$, rel $\partial$, so that it maps each $D_2 \times t$ isometrically to itself. We can now assume that $\widetilde{f}$ maps each $S^2 \times t$ to itself isometrically. In particular, it always maps a pair of antipodal points of $S^2 \times t$ to another pair of antipodal points.

Let $\widetilde{f}_t : S^2 \to S^2$ be the restriction of $\widetilde{f}$ on the level sphere $S^2 \times t$, i.e, $\widetilde{f}_t$ satisfies

$$\widetilde{f}(u,t) = (\widetilde{f}_t(u),t).$$

Define an isotopy $g_s$ on $Y = S^2 \times [0, 1]$ by

$$g_s(u, t) = \begin{cases} (u, t) & \text{if } t \leq s \\[2mm] (\widetilde{f}_{t-s}(u), t) & \text{if } t \geq s \end{cases}$$

Then $g_s$ is an isotopy from $\widetilde{f}$ to the identity map, rel $S^2 \times 0$. Moreover, the restriction of $g_s$ to $S^2 \times 1$ is always an isometry, hence maps antipodal pairs to antipodal pairs. Therefore, $g_s$ induces an isotopy $h_s$ on the quotient space $X = Y/\eta$, deforming $f$ to the identity map rel $\partial X$.  $\square$

Recall from [Sc, Page 446] that there is a unique orientable Euclidean manifold which admits a Seifert fibration with orbifold $P(2, 2)$. We denote this manifold by $M_{P(2,2)}$. The following theorem is sharp: Each of the four manifolds listed admit some map which do not preserve any Seifert fibration of $M$.

**Theorem 3.11.** *Suppose $M$ is a compact orientable Seifert fiber space which is not $T^3$, $M_{P(2,2)}$, $S^1 \times D^2$, or $T \times I$. Then there is a Seifert fibration $p : M \to X(M)$, so that any orientation preserving homeomorphism on $M$ is isotopic to a fiber preserving homeomorphism with respect to this fibration. Moreover, if $M$ is not a lens space or $S^3$, then the result is true for all homeomorphisms on $M$.*

*Proof.* By the remark after Theorem 3.1, we need only consider manifolds covered by $T^3$, $S^3$ or $S^2 \times \mathbb{R}$. The only manifolds covered by $S^2 \times \mathbb{R}$ are $S^2 \times S^1$ and $RP^3 \# RP^3$, in which case the result follows from Lemma 3.10. So we assume that $M$ is a covered by $T^3$ or $S^3$, i.e. it is a Euclidean or spherical manifold.

If $M = S^3$, all orientation preserving maps are isotopic to the identity, which preserves any Seifert fibration of $S^3$. Lemma 3.4 covers lens spaces. The remaining manifolds have Seifert fibrations $p : M \to X(M)$ with orbifolds $X(M) = S(p, q, r)$, or $S(2, 2, 2, 2)$, or $P(2, 2)$, or $T^2$. The case $X(M) = S(2, 2, 2, 2)$ is proved in Lemma 3.9. The last two cases are excluded from the assumption of the theorem. In the remaining cases, $X(M) = S(p, q, r)$. Since $M$ is Euclidean or spherical, we have

$$\frac{1}{p} + \frac{1}{q} + \frac{1}{r} \geq 1.$$

If $(p, q, r) = (2, 2, k)$, or $(2, 3, 4)$, or $(2, 4, 4)$, this is covered by lemmas 3.5, 3.6 and 3.9, respectively. In the remaining cases, either $(p, q, r) = (3, 3, k)$ for some $k \geq 2$, or $(p, q, r) = (2, 3, k)$ for $k \geq 5$, which have been covered by Theorem 3.2.  $\square$

According to [Ja, Theorem VI.17], the following manifolds admit different Seifert fibrations up to isomorphism.

(a) Lens spaces (including $S^3$ and $S^2 \times S^1$).

(b) Prism manifolds, i.e. Seifert fiber spaces with orbifolds $X(M) = S(2, 2, p)$.

(c) Double of twisted $I$-bundle over the Klein bottle, i.e. a Euclidean Seifert fiber space $M$ with $X(M) = S(2, 2, 2, 2)$.

(d) The solid torus.

(e) The twisted $I$-bundle over the Klein bottle.

The following result supplements Theorems 3.1 and 3.2.

**Corollary 3.12.** *Let $M$ be a Seifert fiber space. Then the Seifert fibration of $M$ is unique up to isotopy if and only if $M$ is not $T^3$, $M_{P(2,2)}$, $T \times I$, or one of the manifolds listed above.*

*Proof.* Suppose $M$ is not one of the manifolds listed above. Then by the above theorem of Jaco, Seifert fibrations of $M$ is unique up to isomorphism, that is, if $p_i : M \to X_i$ are Seifert fibrations for $i = 1, 2$, then $X_1 = X_2$, and there is a homeomorphism $f : M \to M$, such that $p_2 \circ f = p_1$. If $M$ is not $T^3$, $T \times I$ or $M_{P(2,2)}$ either, then by Theorem 3.11, up to isotopy we may assume that the map $f$ is a fiber preserving map with respect to the fibration $p_1$. It follows that the two fibrations are the same.

Conversely, if $M$ is one of the manifolds listed above, then it has non isomorphic Seifert fibrations, so of course they are not isotopic. If $M$ is $T^3$ or $T \times I$, then any simple closed curve on $T \times 0$ is a regular fiber of some Seifert fibration of $M$. Since these curves represent different elements in $H_1(M)$, the fibrations are different, so $M$ has infinitely many non isotopic fibrations. Finally for the manifold $M_{P(2,2)}$, it will be shown in Corollary 4.6 that it has exactly three non isotopic Seifert fibrations.   $\square$

## §4. Maps on $T^3$ and $M_{P(2,2)}$

There are two closed manifolds which are excluded in the statement of Theorem 3.11. One is $T^3$, the other one is the manifold $M_{P(2,2)}$, which is the Euclidean manifold with a Seifert fibration over the orbifold $P(2, 2)$. Maps on $T^3$ are clear. We identify $T^3$ with $\mathbb{R}^3/\mathbb{Z}^3$, and call the induced Euclidean metric on $T^3$ the standard metric. A map $f$ on $T^3$ is a linear map if one of its liftings on $\mathbb{R}^3$ is a linear map of the Euclidean space.

**Lemma 4.1.**  *With respect to the standard metric on $T^3$, every map $f : T^3 \to T^3$ is isotopic to a linear map.*

*Proof.* The fundamental group $\pi_1(T^3) = \mathbb{Z}^3$, which is naturally identified with the $\mathbb{Z}^3$ in $\mathbb{R}^3$. The map $f : T^3 \to T^3$ induces a isomorphism $f_* : \mathbb{Z}^3 \to \mathbb{Z}^3$, which extends to a unique linear map $\widetilde{g} : \mathbb{R}^3 \to \mathbb{R}^3$, inducing a linear map $g : T^3 \to T^3$. Since $g_* = f_*$, and $T^3$ is a Haken manifold, $g$ is isotopic to $f$. $\square$

The main objective of this section is to study maps on $M_{P(2,2)}$. Theorem 4.5 will show that all maps on this manifolds are isotopic to a Euclidean isometry.

Up to isotopy there are exactly two orientation preserving essential simple closed curves $a$ and $b$ on a Klein bottle $K$, where $a$ is a nonseparating curve cutting $K$ into an annulus, and $b$ is a separating curve cutting $K$ into two Möbius bands. Let $N$ be a twisted $I$-bundle over $K$. Since any essential annulus in an $I$-bundle over a surface is isotopic to a union of $I$-fibers, we have

**Lemma 4.2.**  *Up to isotopy there are exactly two essential annuli $A, B$ in $N$, which are the union of $I$-fibers over the above curves $a$ and $b$, respectively.* $\square$

Note that a boundary curve $m$ of $A$ intersects a boundary curve $l$ of $B$ at one point. We call $m$ and $l$ the meridian and longitude of the torus $\partial N$. Let $N_1, N_2$ be the two copies of $N$ above. Denote by $A_i, B_i$ the annuli $A, B$ in $N_i$. Similarly for $m_i, l_i$. Gluing $N_1$ and $N_2$ together along their boundary so that $m_i$ is glued to $l_j$, $i \neq j$, we get a manifold $M$. Each of $m_1$ and $l_1$ is a fiber of a fibration of both $N_i$, so we get two Seifert fibrations $p_k : M \to X_k$. The orbifold $p_1(N_1) = U$ is a Möbius band, and $p_1(N_2) = D(2,2)$. Similarly, $p_2(N_1) = D(2,2)$, and $p_2(N_2) = U$. In both cases the orbifolds $X_1 = X_2 = P(2,2)$. Since there are only two fibrations on $N_i$ [WW, Lemma 1.1], these are the only fibrations in which the torus $\partial N_i$ is a vertical torus, i.e. a union of fibers. However, one can show that there is yet another Seifert fibration of $M$ over $P(2,2)$, in which the torus $\partial N_i$ is a horizontal torus, i.e. it intersects each fiber transversely.

Let $A^1 = A'_1 \cup A''_1$ be two parallel copies of $A_1$, let $B^2 = B'_2 \cup B''_2$ be two parallel copies of $B_2$. Since the boundary slope of $B^2$ is $l_2$, which is the same as $m_1$, the same as the boundary slope of $A^1$, we can choose $B^2$ so that $\partial B^2 = \partial A^1$. One can see that $T_1 = A^1 \cup B^2$ is a torus, which is a vertical torus in the fibration $p_1$, and is a horizontal torus in the fibration $p_2$ (i.e. it is transverse to the fibers). Similarly we can define a torus $T_2 = A^2 \cup B^1$, which is horizontal in $p_1$ and vertical in $p_2$.

Note that since the fibration

$$p_i : M \to P(2,2)$$

has a horizontal torus, the Euler number of the fibration is 0, hence $M$ is the Euclidean manifold $M_{P(2,2)}$. Let $T_3 = \partial N_i$. It is a vertical torus in both fibrations $p_i$. We now have three tori $T_1$, $T_2$ and $T_3$. The image of $T_1$ and $T_3$ on the orbifold $P(2,2) = U \cup D(2,2)$ are shown in Figure 4.1, where the image of $T_3$ is the common boundary of $U$ and $D(2,2)$, and the image of $T_1$ is the union of the four arcs. $T_2$ is a horizontal torus, which double covers $P(2,2)$. We will see below that there is a Euclidean metric on $M$ such that the $T_i$'s are mutually perpendicular totally geodesic surfaces.

Figure 4.1

**Lemma 4.3.** *The tori $T_1, T_2, T_3$ are mutually non isotopic. Any essential torus $T$ in $M$ is isotopic to one of the $T_i$.*

*Proof.* The tori intersect each other in a nonempty set of essential curves, so by a Lemma of Waldhausen [Wa, Prop. 5.4], if $T_i$ and $T_j$ are isotopic, then they contain subsurfaces which are parallel in $M$. Cutting along $T_i \cup T_j$, one can check that no such subsurfaces exist.

To prove the second statement, assume $T$ is an essential torus which is not isotopic to $T_3$. By an isotopy we may assume that $T \cap T_3$ are essential curves on $T$ and $T_3$, and each component of $T \cap N_i$ is an essential annulus in $N_i$. By Lemma 4.2, each component of $T \cap N_i$ must be isotopic to $A_i$ or $B_i$, so one can isotop $T$ to a vertical torus in one of the fibrations, say $p_1$. Thus $p_1(T)$ is an orientation preserving essential curve $C$ on the orbifold $X_1 = P(2,2) = U \cup D(2,2)$, and each component of $C \cap U$ and $C \cap D(2,2)$ is an essential arc. Put $k = |C \cap U|$. One can show that if $k = 1$ then $T$ would be a Klein bottle, which is absurd, and if $k \geq 3$ then the curve $C$ is not connected, which contradicts the fact that $T$ is a torus. Thus both $C \cap U$ and $C \cap D(2,2)$ have exactly two components, and the curve $C$ is isotopic to $p_1(T_1)$ on $X_1$. Hence in this case $T$ is isotopic to $T_1$.   $\square$

**Lemma 4.4.** *Let* $f : M \to M$ *be a map such that* $f(T_i)$ *is isotopic to* $T_i$ *for* $i = 1, 2, 3$. *Then* $f$ *is isotopic to a map* $g$ *such that* $g(T_i) = T_i$ *for all* $i$.

*Proof.* By an isotopy we may assume that $f(T_3) = T_3$. Note that $A^1 = T_1 \cap N_1$ consists of two annuli parallel to the annulus $A_1$, which is nonseparating, so if $f$ switches the two sides of $T_3$, then $f(T_1) \cap N_2 = f(T_1 \cap N_1)$ consists of annuli parallel to $A_2$, which is, up to isotopy, the only nonseparating annuli in $N_2$. By the construction of $T_i$ we see that $f(T_1)$ is isotopic to $T_2$, which by Lemma 4.3 is not isotopic to $T_1$, a contradiction. Therefore $f$ must map each $N_i$ to itself. Since up to isotopy $A_1$ is the only essential nonseparating annulus in $N_1$, we can isotop $f$ so that $f(A_1) = A_1$. In particular $f$ preserves the fiber $\partial A_1$ of the fibration $p_1$, so it can be isotoped to a fiber preserving map with respect to the fibration $p_1$. Recall that the orbifold $X_1$ of the fibration $p_1$ is a union of $U$ and $D(2,2)$ over the boundary curve $C$, which is the image of $T_3$. The orbifold map $\widehat{f} : X_1 \to X_1$ maps the two sides of $C$ to itself. It is easy to show that such a map on $X_1$ is isotopic to an identity map. Therefore we may isotop $f$ so that it induces the identity map on $X_1$. Since $T_1$ are $T_3$ are vertical tori, we have $f(T_i) = T_i$ for $i = 1, 3$. It remains to show that we can isotop $f$ via a fiber-wise isotopy so that $f(T_2) = T_2$.

Recall that $T_2 = A^2 \cup B^1$, where $A^2 = A_2' \cup A_2''$ is a pair of annuli parallel to the nonseparating annulus $A_2$ in $N_2$, and $B^1 = B_1' \cup B_1''$ is a pair of horizontal annuli parallel to the separating annulus $B_1$ in $N_1$. We first show that $f$ can be isotoped, via a fiber-wise isotopy, so that $f(B^1) = B^1$.

Cutting $N_1$ along $B^1$, we get three components, each of which is an $I$-bundle over annulus or Möbius band. Since $f$ is a fiber preserving map, $f(B^1)$ is also a horizontal surface, which was cut into horizontal pieces by $B^1$. One can then isotop $f(B^1)$ via a fiber-wise isotopy to reduce the number of components in $f(B^1) \cap B^1$. At the end, $f(B^1)$ will be disjoint from $B^1$, so a further fiber-wise isotopy will move

$f(B^1)$ to $B^1$.

The above fiber-wise isotopy on $N_1$ can be extended to a fiber-wise isotopy of $M$. It remains to show that there is a fiber-wise isotopy on $N_2$, which is the identity on $\partial N_2$, deforming $f(A^2)$ to $A^2$. There are two vertical Möbius bands $U_1$, $U_2$ in $N_2$, cutting $N_2$ into a fibered solid torus. The intersection of $A^2$ with $U_i$ are essential arcs $\alpha$ in $U_i$, similarly for $\alpha' = f(A^2) \cap U_i$. Moreover, since $f(\partial A^2) = \partial A^2$, we have $\partial \alpha = \partial \alpha'$. The key fact here is that any two essential arcs in a Möbius band with the same boundary are rel $\partial$ isotopic to each other. Since both $\alpha$ and $\alpha'$ are horizontal with respect to the induced fibration on $U_i$, this isotopy can be made to be fiber-wise, and extended to a fiber-wise isotopy of $N_2$ rel $\partial N_2$. We can thus isotop $f$ via a fiber-wise isotopy so that $f(\alpha) = \alpha$. After cutting along $U_i$, the surfaces $A^2$ and $f(A^2)$ become two sets of horizontal disks in a fibered solid torus $V$, and they have the same boundary. Therefore we can further fiber-wise isotop $f$ in $V$ rel $\partial V$ to a map $g$ such that $g(A_2) = A_2$.   □

Denote by $I_n$ the interval $[0, n]$. We now identify $N_i$ with the quotient space

$$N_i = I_1 \times I_{\frac{1}{2}} \times I_{\frac{1}{4}} / \sim,$$

where the identification $\sim$ is given by

$$(0, y, z) \sim (1, y, z) \qquad \text{and} \qquad (x, 0, z) \sim (1 - x, \frac{1}{2}, \frac{1}{4} - z).$$

Let $A^1$ be the two annuli at level $x = 0$ and $x = \frac{1}{4}$, and let $B^1$ be the annuli corresponding to the four disks at levels $y = \frac{1}{8}, \frac{3}{8}, \frac{5}{8}, \frac{7}{8}$. Similarly for $A^2$ and $B^2$ in $N_2$. Note that they becomes two annuli after the identification. The boundaries of $A^1$ and $B^1$ cuts $\partial N_1$ into squares with edge length $\frac{1}{4}$, so we can glue $N_1$ to $N_2$ via an isometry on the boundary, which sends $\partial A_i$ to $\partial B_j$, $i \neq j$. The product metric on $I_1 \times I_{\frac{1}{2}} \times I_{\frac{1}{4}}$ now induces a Euclidean metric on $M$, in which each $T_i$ is a totally geodesic torus in $M$, $T_i$ and $T_j$ intersect perpendicularly, and $T = T_1 \cup T_2 \cup T_3$ gives a cell decomposition of $M$, with all 1-cells having length $\frac{1}{4}$, all 2-cells regular squares, and all 3-cells regular cubes.

**Theorem 4.5.** *Any map $f$ on the manifold $M = M_{P(2,2)}$ is isotopic to an isometry with respect to the above metric.*

*Proof.* The map $f : M \to M$ sends each $T_i$ to a torus which, by Lemma 4.3, is isotopic to one of the $T_j$. Therefore $f$ induces a permutation $\tau_f$ on the isotopy classes of $T_1, T_2, T_3$. Let $\mathcal{H}(M)$ be the homeotopy group of $M$, i.e. the group of all isotopy classes of homeomorphisms of $M$. Then we have a map $\tau : \mathcal{H}(M) \to S_3$, where $S_3$ is the permutation group of 3 elements. Let $\mathcal{I}(M)$ be the group of isometries of $M$, and let $\mu : \mathcal{I}(M) \to \mathcal{H}(M)$ be the map induced by inclusion. There is an obvious isometry $\rho_{12}$ switching the two $N_i$, which maps $T_3$ to itself, and switches $T_1$ and $T_2$. The torus $T_1$ also cuts $M$ into two twisted $I$-bundle $N_1'$ and $N_2'$. By examining the above metric, one can see that there is an isometry $\rho_{13}$ which maps $N_i'$ to $N_i$, so it fixes $T_2$ and switches $T_1$ and $T_3$. It follows that the composition $\tau \circ \mu : \mathcal{I}(M) \to S_3$ is a surjective map. To prove the theorem, it remains to show that any isotopy class $[f]$ in $ker(\tau)$ is in the image of $\mu$. Now $[f] \in ker(\tau)$ means that $f(T_i)$ is isotopic to $T_i$ for each $i$. By Lemma 4.4 we can isotop $f$ so that $f(T_i) = T_i$ for all $i$. Thus $f$ preserves the cell structure coming from $\cup T_i$. Since all the cells of the same dimension are isometric, we can isotop $f$ to an isometry, first on the 1-cells, then on the 2-cells, and finally on the 3-cells. The map is then an isometry of $M$.   $\square$

**Corollary 4.6.** *The manifold $M_{P(2,2)}$ has exactly three Seifert fibrations up to isotopy.*

*Proof.* We have seen that, up to isotopy, each Seifert fibration has two vertical essential tori, and each essential torus is a vertical torus for two different Seifert fibrations. Since there are exactly three essential tori in $M$, the result follows.   $\square$

§5. Fixed points and standard maps

Suppose $B^n$ is the unit ball in $\mathbb{R}^n$, and $f : B^n \to \mathbb{R}^n$ is a map with the origin $x$ as the only fixed point. Then the fixed point index of $x$ is defined as

$$ind(f, x) = \text{degree}(\varphi),$$

where the map

$$\varphi : S^{n-1} = \partial B^n \to S^{n-1}$$

is defined by

$$\varphi(y) = \frac{y - f(y)}{||y - f(y)||}.$$

Since a manifold is locally Euclidean, this defines $ind(f, x)$ for isolated fixed point of maps $f : M \to M$ on closed manifolds.

Suppose $f$ is a map on a closed manifold $M$, and $C$ is a union of some components of $Fix(f)$. Then $C$ is compact, and there is an open neighborhood $V$ of $C$ containing no other fixed point of $f$. Perturb $f$ inside $V$ to a map $g$ so that it has only isolated fixed points in $V$. Then $ind(f, C)$ can be defined as the sum of $ind(g, x)$ over all fixed points of $g$ in $V$. It can be shown that this is independent of the perturbation.

Two fixed points $x, y$ of a map $f : M \rightarrow M$ are equivalent if there is a path $\gamma$ from $x$ to $y$, such that $f \circ \gamma$ and $\gamma$ are rel $\partial$ homotopic. In particular, all points in a component of $Fix(f)$ are equivalent. An equivalent class $C$ of fixed points of $f$ is called a Nielsen class. A Nielsen class $C$ is essential if $ind(f, C) \neq 0$; otherwise it is inessential. The Nielsen number $N(f)$ is defined as the number of essential fixed point classes in $Fix(f)$. The crucial property of $N(f)$ is that it is a homotopy invariant. In particular, it gives a lower bound for the number of fixed points for all maps homotopic to $f$.

Figure 5.1

**Example 5.1.** A local map near an isolated fixed point $x$ of $f$ is of *pseudo-Anosov type* if it is the same as the restriction of a pseudo-Anosov map near a fixed point. Figure 5.1 shows the picture of some pseudo-Anosov type local maps. In Figure 5.1(a), the map is a linear map, expanding along the solid line direction, and shrinking along the dotted line direction. One can compute that $ind(f, x) = -1$. More generally, if $f$ is a pseudo-Anosov type map near $x$, as shown in Figure 5.1(b), then $ind(f, x) = 1 - p$, where $p$ is the number of prongs of the singular foliation at $x$. If $f$ is the composition of the above map with a rotation sending a prong to another, the index is 1. In Figure 5.1(c), (d) and (e), the maps are orientation reversing. These are the compositions of maps shown in Figure 5.1(b) and a reflection along an invariant line $l$, which is the union of two singular leaves of the invariant foliations of $f$. The index is 0 if the prong number $p$ is odd (Figure 5.1(c)), $-1$ if $p$ is even and the invariant line is a solid line (Figure 5.1(d)), and 1 otherwise (Figure 5.1(e)). See [JG, Lemma 3.6] for more details. Similarly, a local map in a neighborhood of a boundary point is of pseudo-Anosov type if it is the restriction of a global pseudo-Anosov map in a neighborhood of a boundary fixed point.

The combination of the above example and the following product formula will be useful in computing the index of fixed points for maps of 3-manifolds. See [J1, Lemma I.3.5(iv)] for a proof.

**Lemma 5.2.** *Suppose that $U_i$ is an open neighborhood of $x_i$ in a compact topological space $X_i$, and $f_i : U_i \to X_i$ is a map with $x_i$ as the only fixed point. Let $x = (x_1, x_2)$, and let $f$ be the product map $f = f_1 \times f_2 : X_1 \times X_2 \to X_1 \times X_2$. Then*

$$ind(f, x) = ind(f_1, x_1) \times ind(f_2, x_2). \quad \square$$

If $A$, $B$ are subsets of $M$, then $\gamma : I \to M$ is a path from $A$ to $B$ if $\gamma(0) \in A$, and $\gamma(1) \in B$. If $\gamma_0, \gamma_1$ are paths from $A$ to $B$, and there is a homotopy $\varphi_t : I \to M$ from $\gamma_0$ to $\gamma_1$ such that $\varphi_t(0) \in A$ and $\varphi_t(1) \in B$ for all $t$, then we say that $\gamma_0$ is $(A, B)$ *homotopic to* $\gamma_1$, and write

$$\gamma_0 \sim \gamma_1 \qquad rel\ (A, B).$$

When $A, B$ are one point sets, this is also written as

$$\gamma_0 \sim \gamma_1 \qquad rel\ \partial.$$

A set $B$ is *f-invariant* if $f(B) = B$. Suppose $f : M \to M$ is a map, and $A, B$ are $f$-invariant sets of $M$. If there is a path $\gamma$ from $A$ to $B$ such that $\gamma \sim f \circ \gamma$ rel $(A, B)$, then we say that $A, B$ are *f-related*, denoted by

$$A \sim_f B.$$

In particular, two fixed points $x$ and $y$ of a map $f$ are in the same fixed point class if and only if they are $f$-related.

A covering $p : \widetilde{M} \to M$ is *characteristic* if the corresponding subgroup $p_* \pi_1 \widetilde{M}$ is a characteristic subgroup of $\pi_1 M$, i.e, it is invariant under any automorphism of $\pi_1 M$. Recall that if $G$ is a finite index subgroup of a finitely generated group $H$, then there is a finite index characteristic subgroup $G'$ of $H$, such that $G' \subset G$ (see [???]). It follows that given any finite covering $p : \widetilde{M} \to M$ of compact 3-manifold, there is a finite covering $q : \widehat{M} \to \widetilde{M}$ so that $p \circ q : \widehat{M} \to M$ is a characteristic covering. Notice that if $p : \widetilde{M} \to M$ is a characteristic covering, then for any homeomorphism $f : M \to M$ with $y = f(x)$, and any points $\widetilde{x}$ and $\widetilde{y}$ covering $x$ and $y$ respectively, there is a lifting $\widetilde{f}$ of $f$ such that $\widetilde{f}(\widetilde{x}) = \widetilde{y}$.

**Lemma 5.3.** *Suppose that $A$ and $B$ are $f$-invariant subsets of $M$ which are $f$-related by a path $\gamma$. Let $p : \widetilde{M} \to M$ be a characteristic covering map, and $\widetilde{A}$ a component of $p^{-1}(A)$. Then there is a lifting $\widetilde{f}$ of $f$, and a lifting $\widetilde{\gamma}$ of $\gamma$, which $\widetilde{f}$-relates $\widetilde{A}$ to a component $\widetilde{B}$ of $p^{-1}(B)$.*

*Proof.* First choose $\widetilde{\gamma}$ that covers $\gamma$. The homotopy $h_t : \gamma \sim f \circ \gamma$ lifts to a homotopy $\widetilde{h}_t$ on $\widetilde{M}$, with $\widetilde{h}_0(0) = \widetilde{\gamma}(0)$. Then $p \circ \widetilde{h}_1(0) = h_1(0) = f \circ \gamma(0)$, so $\widetilde{h}_1(0) \in p^{-1}(f \circ \gamma(0))$. Since $p$ is characteristic, there is a lifting $\widetilde{f}$ of $f$ such that $\widetilde{f}(\widetilde{\gamma}(0)) = \widetilde{h}_1(0)$. Thus $\widetilde{h}_t$ is a homotopy $\widetilde{\gamma} \sim \widetilde{f} \circ \widetilde{\gamma}$ rel $(\widetilde{A}, \widetilde{B})$, where $\widetilde{B}$ is the component of $p^{-1}(B)$ containing $\widetilde{\gamma}(1)$.   $\square$

**Definition 5.4.** A *normal structure* on a properly embedded arc or circle $C$ in a 3-manifold $M$ is an identification of a neighborhood $N(C)$ with $C \times D$, where $D$ is a disk, such that $C$ is identified with $C \times 0$. Suppose $C$ is a proper 1-dimensional component of $Fix(f)$. We say that $f$ *preserves a normal structure* on $C$ if there is a normal structure $N(C) = C \times D$ on $C$, such that $f$ restricting to a smaller neighborhood $N_1(C) = C \times D_1$ has the property that $f(x \times D_1) \subset x \times D$ for any $x \in C$.

**Lemma 5.5.** *Suppose $f : M \to M$ is an orientation preserving map. If $C$ is a circle component of $Fix(f)$, and $f$ preserves the normal structure on $C$, then $f$ can be isotoped in a neighborhood of $C$ to a map $g$ such that $Fix(g) = Fix(f) - C$.*

*Proof.* Let $N_1(C) = C \times D_1$ be as in Definition 5.4. Identify $f(N_1(C))$ with $S^1 \times D'$, where $D'$ is the unit disc. Define an isotopy of $M$ by $\varphi_t(x, p) = (x + t(1 - ||p||), p)$ on $f(N_1(C))$ and identity outside, where $||p||$ is the norm of $p$. Since $f$ maps each $x \times D'$ into $x \times D$, the composition $g = \varphi_1 \circ f$ has no fixed points in $N_1(C)$.    $\square$

**Definition 5.6.** A map $f : M \to M$ is said to be of *flipped pseudo-Anosov type* at an isolated fixed point $x$ if there is a neighborhood $N(x) = D \times I$ of $x$, such that the restriction of $f$ on a smaller neighborhood $D_1 \times I$ is a product $f_1 \times f_2$, where $f_1$ is a pseudo-Anosov type map given in Example 5.1, and $f_2$ is the involution $f_2(t) = 1 - t$.

Notice that in this case, if we write $x = (x_1, x_2)$, then by Lemma 5.2 we have $ind(f, x) = ind(f_1, x_1) \cdot ind(f_2, x_2) = ind(f_1, x_1)$, because in this case $ind(f_2, x_2) = 1$.

**Lemma 5.7.** *Suppose $f : M \to M$ is an orientation preserving map. If $f$ is of flipped pseudo-Anosov type at $x \in Fix(f)$, and $ind(f, x) = 0$, then $f$ can be isotoped in a neighborhood of $x$ to a map $g$ such that $Fix(g) = Fix(f) - x$.*

*Proof.* Let $f = f_1 \times f_2 : D_1 \times I \to D \times I$ be as in the definition, and let $x = (x_1, x_2)$. Then $ind(f_1, x_1) = ind(f, x) = 0$. By [JG, §2.2, Remark 2], we can isotop $f_1$ rel $\partial D_1$ to a map $g_1$ which has no fixed point. Extend this isotopy to an isotopy of $f$ in $D_1 \times I$, rel $\partial(D_1 \times I)$. Since $f$ exchanges the two sides of $D_1 \times 0$, this will not create new fixed points.

**Definition 5.8.** Suppose $M$ is a compact 3-manifold with torus boundaries. A map $f : M \to M$ is *standard on boundary* if for any component $T$ of $\partial M$, the map $f|_T$ is of one of the following types.

(i) A fixed point free map;

(ii) A periodic map with isolated fixed points;

(iii) A fiber preserving, fiber orientation reversing map with respect to some $S^1$ fibration of $T$.

**Definition 5.9.** (1) A map $f$ on a compact 3-manifold $M$ is said to have *FR-property* (Fixed-point Relating Property) if the following is true:

If $A \in Fix(f)$ and $B$ is either a fixed point of $f$ or an $f$-invariant component of $\partial F$, and $A, B$ are $f$-related by a path $\gamma$, then $\gamma$ is $(A, B)$ homotopic to a path in $Fix(f)$.

(2) $f$ is a *Type 1 standard map* if (i) $f$ has FR-property, (ii) $f$ is standard on boundary, (iii) $Fix(f)$ consists of isolated points, and (iv) $f$ is of flipped pseudo-Anosov type at each fixed point.

(3) $f$ is a *Type 2 standard map* if it satisfies (i), (ii) above, as well as (iii) $Fix(f)$ is a properly embedded 1-dimensional submanifold, and (iv) $f$ preserves a normal structure on $Fix(f)$.

**Lemma 5.10.** *Suppose $M$ is a closed 3-manifold. If $f : M \to M$ is a type 1 standard map, then $f$ is isotopic to a map $g$ with exactly $N(f)$ fixed points. If $f$ is a type 2 standard map, then it is isotopic to a fixed point free map.*

*Proof.* The FR-property implies that each fixed point class of $f$ is path connected. In particular, if $f$ is a Type 1 standard map on a closed manifold, then each fixed point is a Nielsen class, and by Lemma 5.7 we can isotop $f$ in a neighborhood of inessential fixed points to eliminate them, deforming $f$ to a map $g$ with exactly $N(f)$ fixed points. If $f$ is a type 1 standard map, since $M$ is closed, all components of $Fix(f)$ are circles, which can be eliminated using Lemma 5.5.  $\square$

The following result about surface maps is due to Jiang and Guo [JG, Lemma 1.2 and Lemma 2.2] when $F$ is hyperbolic. The result is well known if $F$ is Euclidean, and can be proved in the same way.

**Lemma 5.11.** *(1) If $\varphi : F \to F$ is a periodic map on a hyperbolic or Euclidean surface, then it has FR-property.*

*(2) If $\varphi : F \to F$ is an Anosov or pseudo-Anosov map, then each fixed point of $\varphi$ is an isolated point, and it has FR-property, except that two fixed points on the same boundary component $C$ of $F$ are $\varphi$-related if $\varphi$ is isotopic to id on $C$. In particular, if $\varphi$ is orientation reversing, then it has FR-property.*  $\square$

The following simple fact is useful in detecting equivalent fixed points.

**Lemma 5.12.** *Let $f$ be a map on $M$, and let $p : \widetilde{M} \to M$ be the universal covering of $M$. Then two points $x_1, x_2 \in Fix(f)$ are equivalent if and only if there is a lifting $\widetilde{f} : \widetilde{M} \to \widetilde{M}$ of $f$, and points $\widetilde{x}_i \in Fix(\widetilde{f})$, so that $p(\widetilde{x}_i) = x_i$. In particular, if $Fix(\widetilde{f})$ is connected for all liftings $\widetilde{f}$ of $f$, then no two components of $Fix(f)$ are equivalent.*

*Proof.* The necessity is a special case of Lemma 5.3. The sufficiency follows by projecting a homotopy $\gamma' \sim \widetilde{f} \circ \gamma'$ to one on $M$, where $\gamma'$ is a path from $\widetilde{x}_1$ to $\widetilde{x}_2$. Such a homotopy exists because $\widetilde{M}$ is simply connected.   $\square$

**Lemma 5.13.** *Let $f : \mathbb{R}^3 \to \mathbb{R}^3$ be a linear map defined by a matrix $A$, and suppose that $f$ has a single fixed point at the origin $x$. Then $ind(f, x) = det(E - A) \neq 0$, where $E$ is the identity matrix.*

*Proof.* This is a special case of [J1, Properties 3.2(2)]. Notice that since $f$ fixes no 1-dimensional subspaces, $det(E - A) \neq 0$.   $\square$

**Lemma 5.14.** *Suppose $f$ is an orientation reversing map on a torus $T$. Then $f$ is isotopic to either a fixed point free map or an Anosov map.*

*Proof.* By an isotopy we may assume that $f$ is either periodic, reducible, or Anosov. If it is periodic and $Fix(f)$ is nonempty, then it is an involution, so it is also reducible. For reducible map, its fixed point set consists of two circles, which can be eliminated by a small flow along the circle direction.   $\square$

§6. Fiber preserving maps

Let $f : M \to M$ be an orientation preserving, fiber preserving map on a Seifert fiber space $M$, let $p : M \to X(M)$ be the projection map, and let $\widehat{f} : X(M) \to X(M)$ be the induced map on the orbifold. In this section we will consider the case that $\widehat{f} : X(M) \to X(M)$ is either periodic, Anosov, or pseudo-Anosov.

Consider the restriction of $f$ on $p^{-1}(W)$ for a component $W$ of $Fix(\widehat{f})$. For each $x \in W$, the fiber $L_x = p^{-1}(x)$ is $f$-invariant. Since $W$ is connected, it is clear that $f$ preserves the orientation of a fiber over $W$ if and only if it preserves the orientation of all fibers over $W$. In this case we say that $f$ is fiber orientation preserving over $W$. Otherwise it is fiber orientation reversing over $W$.

**Lemma 6.1.** *Suppose $f$ is fiber orientation preserving over some component $W$ of $Fix(\widehat{f})$. Let $M' = p^{-1}(N(W))$. Then $f$ can be isotoped, via a fiber preserving isotopy rel $M - M'$, to a map which has no fixed point in $M'$.*

*Proof.* By a small fiber preserving isotopy rel $M - M'$ we may modify $f$ to a map $g$ such that the orbifold map $\widehat{g}$ has only isolated fixed points in $N(W)$. Thus $g$ maps only finitely many fibers $C_1, \ldots C_n$ of $M'$ to themselves.

For each $C_i$, since $g$ is fiber orientation preserving, $g|_{C_i}$ can be isotoped to a fixed point free map. The isotopy can be extended to a fiber-wise isotopy supported in a small neighborhood of $C_i$, eliminating all fixed points in $N(C_i)$. The result follows by performing such an isotopy for each $C_i$.   $\square$

Notice that if $M$ is a connected orientable circle bundle over an orientable surface, then all fibers can be coherently oriented, so a fiber preserving map $f : M \to M$ either preserves or reverses the orientation on all fibers.

**Lemma 6.2.** *Suppose $q : M \to F$ is a connected orientable circle bundle over an orientable compact surface $F$ with $\chi(F) \leq 0$, and suppose $f : M \to M$ is an orientation preserving, fiber orientation reversing map, such that $\widehat{f} : X(M) \to X(M)$ is periodic, Anosov or pseudo-Anosov. Then $f$ has FR-property.*

*Proof.* Consider the case that $A, B$ are fixed points of $f$, and are $f$-related by a path $\gamma$. (The case that $B$ is an $f$-invariant boundary component of $f$ is similar.) Then $\widehat{A} = q(A)$ and $\widehat{B} = q(B)$ are fixed points of $\widehat{f}$ on $X(M)$, and they are $\widehat{f}$-related by $\widehat{\gamma} = q \circ \gamma$. Since $f$ is fiber orientation reversing, $\widehat{f}$ is orientation reversing on $F$. Hence by Lemma 5.11, $\widehat{\gamma}$ is $(\widehat{A}, \widehat{B})$ homotopic to a path $\widehat{\beta}$ in $Fix(\widehat{f})$. As $f$ is fiber orientation reversing, $Fix(\widehat{f})$ is double covered by $Fix(f)$, so there is a unique path $\beta$ in $Fix(f)$ covering $\widehat{\beta}$, with $\beta(0) = \gamma(0)$. The path $\beta^{-1} \cdot \gamma$ projects to the null-homotopic loop $\widehat{\beta}^{-1} \cdot \widehat{\gamma}$ on $F$, so it is rel $\partial$ homotopic to a path $\alpha$ in the fiber $L$ over $\widehat{\gamma}(1)$. Since

$$\alpha \sim \beta^{-1} \cdot \gamma \sim f \circ (\beta^{-1} \cdot \gamma) \sim f \circ \alpha \qquad \text{rel } \partial,$$

and since $f$ reverses the orientation of $L$, this is true only if $\alpha$ is a null-homotopic loop in $L$. It follows that $\gamma \sim \beta$ rel $\partial$.   $\square$

**Proposition 6.3.** *Suppose $f : M \to M$ is a fiber preserving, orientation preserving map on a Seifert fiber space $M$ with hyperbolic or Euclidean orbifold.*

*(1) If $\widehat{f} : X(M) \to X(M)$ is Anosov or pseudo-Anosov, then $f$ is isotopic to a Type 1 standard map.*

*(2) If $\widehat{f} : X(M) \to X(M)$ is periodic, then $f$ is isotopic to a Type 2 standard map.*

*Proof.* By Lemma 6.1 we may assume that $f$ reverses the orientation of any fiber which contains a fixed point of $f$. In particular, if $T$ is a component of $\partial M$, then either $f|_T$ is fixed point free, or it is a fiber preserving, fiber orientation reversing map. Hence $f$ is standard on boundary.

Suppose $A \in Fix(f)$ is $f$-related to $B$ by a path $\gamma$, where $B$ is a fixed point of $f$ or an $f$-invariant component of $\partial M$. By [Sc, Theorem 2.5], there is a finite covering $q : X(M') \to X(M)$ of orbifold such that $X(M')$ is an orientable surface. The pull-back of the Seifert fibration $p : M \to X(M)$ via $q$ gives a 3-manifold $M'$ and a Seifert fibration $q' : M' \to X(M')$. Since $X(M')$ is an orientable surface, the fibration $q'$ is an orientable circle bundle over the orientable surface $F = X(M')$. By the remark before Lemma 5.3, after passing to a further finite covering if necessary, we may assume that $q'$ is characteristic.

Let $f' : M' \to M'$ be a lifting of $f$ which fixes a point $A'$ with $q'(A') = A$. By Lemma 5.3, $A'$ is $f'$-related to some component $B'$ of $(q')^{-1}(B)$, by a path $\gamma'$ which is a lifting of $\gamma$. Since $A$ is a fixed point of $f$, by the above assumption $f$ reverses the orientation of the fiber containing $A$. As a lifting of $f$, the map $f'$ reverses the orientation of the fiber containing $A'$, and hence the orientation of all fibers because the fibers of $M'$ can be coherently oriented. Since the orbifold map $\widehat{f'} : X(M') \to X(M')$ is a lifting of the (periodic, Anosov or pseudo-Anosov) map $\widehat{f}$, it is also periodic, Anosov or pseudo-Anosov. By Lemma 6.2, $f'$ has FR-property. Hence there is a homotopy $h'_t$ deforming $\gamma'$ to $\beta' \subset Fix(f')$ rel $(A', B')$. The projection of $h'_t$ on $M$ is a homotopy $q' \circ h'_t : \gamma \sim q' \circ \beta'$ rel $(A, B)$. Since $q' \circ \beta'$ is a path in $Fix(f)$, this proves that $f$ has FR-property.

Suppose $\widehat{f}$ is Anosov or pseudo-Anosov. Then $Fix(\widehat{f})$, hence $Fix(f)$, are isolated points. For each point $x \in Fix(f)$, choose a product neighborhood $N(x) = D \times I$, with each $p \times I$ in a Seifert fiber of $M$. By a fiber-wise isotopy we may assume that $f$ is a product $f_1 \times f_2$ in a smaller neighborhood $N_1(x) = D_1 \times I$, where $f_2$ is the involution on $I$. The map $f_1 : D \to D$ is a lifting of $\widehat{f}|_{p(D_1)}$, which is of pseudo-Anosov type. Therefore $f$ is a type 1 standard map.

Now suppose $\widehat{f}$ is periodic. Then $Fix(\widehat{f})$ consists of geodesics, so $Fix(f)$ is a 1-dimensional submanifold. Let $C$ be a component of $Fix(f)$. Let $N(C)$ be a regular neighborhood of $C$ which is an $I$-bundle induced by the Seifert fibration, over a surface $A$, which is either an annulus, a Möbius band, or a rectangle. Let $q : N(C) \to A$ be the quotient map pinching each $I$-fiber to a point. The map $f$ induces a map $f_1 : A_1 \to A$, where $A_1$ is a subsurface of $A$ containing $q(C)$. Notice that $p : M \to X(M)$ induces a map $r : A \to X(M)$, which is locally a covering map of orbifold, and $f_1$ covers the orbifold map $\widehat{f}$, i.e, $r \circ f_1 = \widehat{f} \circ r$. Since $\widehat{f}$ is periodic with a 1-dimensional fixed point set, it is an involution, hence $f_1$ is locally an involution near each fixed point. Therefore, we may choose a smaller $A$ if necessary, so that $f_1$ is an involution on $A$. It is now clear that $f_1$ preserves some normal structure of the curve $q(C)$ on $A$, which then lifts to a normal structure on $N(C)$ preserved by $f$. This completes the proof that $f$ is a type 2 standard map.   $\square$

The following result supplements Proposition 6.3. It covers fiber preserving maps on manifolds which have spherical orbifolds.

**Proposition 6.4.** *Suppose that $M$ is a closed Seifert fiber space with orbifold $X(M)$ a sphere with at most three cone points, or a projective plane $P$ with at most two cone points. Let $f : M \to M$ be a fiber preserving map. Then $f$ is isotopic to a fixed point free map.*

*Proof.* By a fiber preserving isotopy we may assume that $\widehat{f} : X(M) \to X(M)$ is periodic. By the argument in the last paragraph of the proof of Proposition 6.3, $f$ is isotopic to a map such that $Fix(f)$ is a 1-manifold, and $f$ preserves a normal structure of $Fix(f)$. By Lemma 5.5 we can isotop $f$ to remove all such fixed points.   $\square$

## §7. Maps on torus bundles

A 1-orbifold $Y$ is an interval with 0, 1 or 2 cone points at its ends, or an $S^1$. In the first case a torus bundle over $Y$ is a $T \times I$, in the second case a twisted $I$-bundle over the Klein bottle, in the third case a union of two twisted $I$-bundles over Klein bottle, and in the fourth case it is a torus bundle over $S^1$. In this section we will study maps on these manifolds with each fiber a torus. In particular, we will prove Theorem 9.1 for maps on Sol manifolds.

**Lemma 7.1.** *If $M$ is a twisted $I$-bundle over Klein bottle, then a map $f : M \to M$ can be isotoped to a type 2 standard map.*

*Proof.* According to [WW, Lemma 1.1], up to isotopy $M$ admits exactly two Seifert fibrations: It is a twisted circle bundle over Möbius band, and it is also a Seifert fiber space over a disk with two singular points of index 2. Since those two Seifert fibrations are not homeomorphic, the homeomorphism $f$ must preserve each Seifert fibration up to isotopy. So we may assume that $f$ is a fiber preserving map. By Lemma 2.3, with a further isotopy we may assume that the orbifold map $\widehat{f}$ is periodic. Hence we can apply Proposition 6.3 to isotop $f$ to a type 2 standard map. Since the restriction of $f$ to the central Klein bottle $K$ is isotopic to a periodic map by Proposition 2.4, the map $f|_{\partial M}$, as a covering of $f|_K$, is also isotopic to a periodic map. By the definition of standard map, $f|_{\partial M}$ is either fixed point free or a periodic map with isolated fixed points.   $\square$

**Lemma 7.2.** *Let $N = T \times I$, where $T$ is a torus. Let $T_i = T \times i$, and let $f_i = f|_{T_i}$. Suppose $f : N \to N$ is a map which is standard on boundary.*

*(1) If $f_i(T_i) \neq T_i$, then $f$ is rel $\partial N$ isotopic to a type 1 standard map $g$. Moreover, if $f_i^2 : T_i \to T_i$ is not isotopic to an Anosov map, then $g$ is fixed point free.*

*(2) If $f_i(T_i) = T_i$, then $f$ is rel $\partial N$ isotopic to a type 2 standard map.*

*Proof.* (1) Suppose $f_i(T_i) \neq T_i$. We may assume that $f$ maps $T = T \times \frac{1}{2}$ to itself. By Lemma 5.14 we may isotop $f' = f|_T$ to a fixed point free or Anosov map. The isotopy can be extended to an isotopy of $f$ rel $\partial N$, modifying $f$ to a map $g$ so that in a neighborhood of $T$ it is the product of $f'$ with an involution. Each point of $Fix(g)$ is an essential fixed point class of $g$, and since $g$ has no invariant boundary component, no point of $Fix(g)$ is $g$-related to a boundary component, so $g$ has FR-property. Therefore $g$ is a type 1 standard map, and $g$ is fixed point free unless $f'$ is Anosov, which is true if and only $f_i^2 : T_i \to T_i$ is isotopic to an Anosov map.

(2) Suppose $f_i(T_i) = T_i$. We separate the argument into several cases.

CASE 1. ($f_i$ has not fixed point.)

In this case we can just isotop $f$ rel $\partial N$ so that for any $x \in (0,1)$, $f$ sends the torus $T \times x$ to $T \times y$ for some $y \in (0, x)$. Then $f$ has no fixed point on $N(T)$, so it is automatically a type 2 standard map.

CASE 2. ($f_i$ is a periodic map with isolated fixed points.)

In this case it follows from [JW, Proposition 3.1] that $f$ is rel $\partial N$ isotopic to an affine map $g$ with respect to some Euclidean structure of $T_times I$, so that each component of $Fix(g)$ is an arc from $T_0$ to $T_1$. (Note that since $f|_{T_i}$ is orientation preserving and is not the identity map, $Fix(f)$ can not be torus or circle on $T_i$, which rules out the other possibilities of [JW, Prop. 3.1].) Since $g$ is affine, it maps each $T \times t$ to itself, hence it preserves a normal structure.

We need to prove that $g$ has FR-property. If $A \in Fix(g)$ is $g$-related to a boundary component $B$ of $N$ by $\gamma$, then $\gamma$ is $(A, B)$ homotopic to the arc in $Fix(g)$ containing $A$. So it remains to show that two points in different components of $Fix(g)$ are not $g$-related. If this were not true, since each component of $Fix(g)$ intersects $T_0$, we can find $x, y \in Fix(g) \cap T_0$, and a path $\gamma$ which $g$-relates $x$ to $y$. The curve $\gamma$ is rel $\partial$ homotopic to a path $\beta$ in $T_0$, so $\beta \sim g \circ \beta$ rel $\partial$. Projecting the homotopy to $T_0$, we see that two fixed points of $g|_{T_0}$ are $g$-related on $T_0$. This is impossible because $g|_{T_0}$ is periodic, and by Lemma 5.11 its fixed points are not $g|_{T_0}$-related.

CASE 3. $f_i$ *is a fiber preserving, fiber orientation reversing map with respect to some $S^1$ fibration of $T_i$.*

First assume that the two fibrations on $T_i$ are isotopic to each other. In this case the fibrations on $T_i$ extends to an $S^1$ bundle structure of $N$ over an annulus $A$, and $f$ is rel $\partial N$ isotopic to a fiber preserving map. The induced map $\widehat{f} : A \to A$ is orientation reversing, so it is rel $\partial A$ isotopic to an involution. The result now follows immediately from Lemma 6.2.

Now assume that the fibers of $T_i$ are not isotopic to each other. Notice that in this case the fibers generates $\pi_1 N$, so $f_i$ is isotopic to an involution. Extend the circle fibration over a neighborhood of $\partial N$, and isotop $f$ rel $\partial N$ to a map $g$ which is fiber preserving on $N(T_i)$, and $g|_{T_i'}$ is an involution, where $T_i' = \partial N(T_i) - T_i$. Let $N'$ be the closure of $N - N(T_1) \cup N(T_2)$. By Case 2 we can further isotop $g$ rel $N(T_1) \cup N(T_2)$ so that $g|_{N'}$ is a type 2 standard map. One can check as in Case 2 that $g$ is a type 2 standard map on $N$. □

Recall that $N(f)$ denotes the Nielsen number of a map $f$.

**Lemma 7.3.** *Suppose $M$ is a torus bundle over $S^1$, and $f : M \to M$ is an orientation preserving map which preserves the torus fibration. Then $f$ can be isotoped to a map $g$ with exactly $N(f)$ fixed points. Furthermore, $N(f) = 0$ unless (i) the induced map $f' : S^1 \to S^1$ is orientation reversing, and (ii) the restriction of $f$ on at least one of the invariant torus fibers is isotopic to an Anosov map.*

*Proof.* If $f' : S^1 \to S^1$ is orientation preserving, then we can deform $f$ via a fiber preserving isotopy, so that $f'$ is fixed point free, hence $f$ is also fixed point free, and the result follows.

If $f'$ is orientation reversing, then $f$ has two invariant torus fibers $T_1$ and $T_2$, and $Fix(f) \subset T_1 \cup T_2$. By the Nielsen Realization Theorem of surface homeomorphisms [JG, Main Theorem], we may isotop $f_i = f|_{T_i}$ so that $N(f_i) = \#Fix(f_i)$ as maps on tori. The isotopy can be extended to a fiber preserving isotopy of $f$, deforming $f$ to a new map $g$ such that the restriction of $g$ on $N(T_i) = T_i \times I$ is a product $f_i \times \varphi$, where $\varphi : I \to I$ is the involution $\varphi(t) = 1 - t$. Since $T_i$ are still the only fibers preserved by $g$, we have $Fix(g) = Fix(f_1) \cup Fix(f_2)$. Let $t = 1/2$ be the fixed point of $\varphi$. Note that since $\varphi$ is an involution, the index of $t$ is 1. Since $g$ is locally a product $f_i \times \varphi$ near each fixed point, by Lemma 5.2 we have

$$ind(f, x) = ind(f_i, x) \times ind(\varphi, t)$$

for any $x \in Fix(f)$. It follows that the index of each fixed point of $g$ is nonzero. We need to show that no two of them are Nielsen equivalent.

Suppose $x, y \in Fix(g)$ are in the the same fixed point class of $g$. Then there is a path $\gamma$ from $x$ to $y$ such that

$$g \circ \gamma \sim \gamma \qquad \text{rel } \partial.$$

Let $\gamma' = q \circ \gamma$. Then $\gamma'$ is a path from $q(x)$ to $q(y)$, and

$$g' \circ \gamma' \sim \gamma' \qquad \text{rel } \partial,$$

where $g'$ is the orbifold map on $S^1$ induced by $g$. Since $g'$ is orientation reversing, this implies that $\gamma'$ is homotopic to a trivial loop. Thus $x$ and $y$ must be in the same $T_i$, and $\gamma$ can be deformed into $T_i$. Now $\pi_1(T_i)$ injects into $\pi_1(M)$, so the homotopy $\gamma \sim f \circ \gamma$ can be deformed into $T_i$, which implies that $x$ and $y$ are in the same fixed point class of $f_i$, contradicting the fact that $\#Fix(f_i) = N(f_i)$. Therefore the number of fixed points of $g$ is exactly $N(f)$.

By Lemma 5.14, $f_i$ is isotopic to a fixed point free or Anosov map. So if neither $f_i$ is isotopic to an Anosov map, then $g$ is fixed point free.   $\square$

Denote by $J$ the 1-orbifold $[-1, 1]$ with two singular points at its ends. If $p : M \to J$ is a torus fibration, then $M$ is a union of two twisted $I$-bundle over Klein bottle, glued together along their boundary.

**Lemma 7.4.** *Suppose $M$ is a torus bundle over the 1-orbifold $J$, and $f : M \to M$ is an orientation preserving map which preserves the torus fibration. Then $f$ is isotopic to a map with no fixed points.*

*Proof.* Let $f' : J \to J$ be the orbifold map induced by $f$. By an isotopy we may assume that $f'$ is either the identity or an involution on $J$. In particular, $f'$ fixes the middle point of $J$, so the middle torus fiber $T$ is invariant under the map $f$.

First assume that $f'$ is an involution on $J$. As in the proof of Lemma 7.3, $f|_T$ is isotopic to either a fixed point free map, or an Anosov map. Now $T$ cuts $M$ into $N_1, N_2$, which are twisted $I$-bundles over Klein bottle. The map $f^2$ sends each $N_i$ to itself, so by Lemma 2.4 it is isotopic to a periodic map. In particular, $f^2|_{\partial N_i} = (f|_T)^2$ is isotopic to a periodic map, so $f|_T$ can not be an Anosov map. Thus $f|_T$ is isotopic to a fixed point free map. The isotopy can be fiber-wise extended over $M$, modifying $f$ to a fixed point free map.

Now assume that $f' : J \to J$ is the identity map. Since $N_i$ is a twisted $I$-bundle over Klein bottle, by Lemma 7.1 we can isotop $f$ so that $f|_{N_i}$ is a Type 2 standard map. By Lemma 7.2 we can isotop $f|_{N(T)}$ rel $\partial N(T)$ to a Type 2 standard map. Thus $Fix(f)$ is a union of 1-manifold in $N_i$ and $N(T)$, with some normal structure preserved by $f$. Since $f|_{N_i}$ and $f|_{N(T)}$ have no isolated fixed points, $Fix(f)$ must be a union of circles in $M$. By Lemma 5.5 one can isotop $f$ to eliminate all such fixed points, and get a fixed point free map.  $\square$

**Proposition 7.5.** *Suppose $M$ is an orientable 3-manifold with Sol geometry and $f : M \to M$ is an orientation preserving map. Then $f$ is isotoped to a fixed point free map.*

*Proof.* A Sol manifold has a structure of torus bundle over 1-orbifold, and by [Sc, Theorem 5.3] it is not a Seifert fiber space. Thus by Lemmas 3.7 and 3.8, we may assume that $f$ is torus fiber preserving. If $M$ is a torus bundle over the 1-orbifold $J$, the result follows from Lemma 7.4. Hence we may assume that $M$ is a torus bundle over $S^1$. By Lemma 7.3, $f$ can be isotoped so that $\#Fix(f) = N(f)$. It remains to show that $N(f) = 0$.

Figure 7.1

Continue with the proof of Lemma 7.3. We have shown that $f$ induces an involution on the orbifold $S^1$, whose two fixed points yield two $f$-invariant tori $T_1, T_2$ of $M$, and $Fix(f) = Fix(f_1) \cup Fix(f_2)$, where $f_i = f|_{T_i}$. We may consider $M$ as a quotient

$$M = (T \times I)/\varphi$$

where

$$\varphi : T \times 0 \to T \times 1$$

is the gluing map. We may assume that $T_1 = T \times 1 = T \times 0$, and $T_2 = T \times \frac{1}{2}$. Thus $f$ pulls back to a torus fiber preserving map

$$\widehat{f} : T \times I \to T \times I$$

which induces an involution on $I$. Put $f' = \widehat{f}|_{T \times 0}$, and $f'' = \widehat{f}|_{T \times 1}$. See Figure 7.1. We have the following commutative diagram.

$$
\begin{array}{ccc}
T \times 0 & \xrightarrow{\ \varphi\ } & T \times 1 \\
{\scriptstyle f'}\big\downarrow & & \big\downarrow{\scriptstyle f''} \\
T \times 1 & \xleftarrow{\ \varphi\ } & T \times 0
\end{array}
$$

Fixing a basis for $H_1(T)$, each map then determines a matrix with respect to this basis. Let $A$ be the matrix of $\varphi$. Since $f'$ and $f''$ are the restriction of $\widehat{f}$, they induce the same matrix $B$. The above commutative diagram implies that

$$BA^{-1} = AB.$$

Since $M$ is a Sol manifold, the gluing map $\varphi$ is an Anosov map, so the matrix $A$ has two *distinct* real eigen values $\lambda_1, \lambda_2$. Let $v_i$ be an eigen vector of $\lambda_i$. With respect to the basis $(v_1, v_2)$ the map $\varphi$ is represented by a diagonal matrix with diagonal entries $\lambda_1$ and $\lambda_2$. Let $B'$ be the matrix representing $f'$ and $f''$ with respect to this new basis. Then

$$B' \begin{pmatrix} \lambda_2 & 0 \\ 0 & \lambda_1 \end{pmatrix} = \begin{pmatrix} \lambda_1 & 0 \\ 0 & \lambda_2 \end{pmatrix} B'$$

Since $\lambda_1 \neq \lambda_2$, by comparing the coefficient in the above equation we have

$$B' = \begin{pmatrix} 0 & b_{12} \\ b_{21} & 0 \end{pmatrix},$$

hence $tr(B) = tr(B') = 0$.

Denote by $h$ the map $f_2 : T_2 \to T_2$. Let

$$h_{i*} : H_i(T_2) \to H_i(T_2)$$

be the induced homomorphism on the $i$-th homology group of $T_2$. The matrix $B$ is also the representation matrix of $h_{1*}$. According to [Ki], the Nielsen number of a torus map is the same as its Lefschetz number:

$$N(h) = tr(h_{0*}) - tr(h_{1*}) + tr(h_{2*}).$$

Clearly $tr(h_{0*}) = 1$. Since $h = f_2$ is orientation reversing on $T_2$, we have $tr(h_{2*}) = -1$. Therefore

$$N(f_2) = N(h) = tr(h_{1*}) = tr(B) = 0.$$

Similarly, one can show that $N(f_1) = 0$. Since $Fix(f) = Fix(f_1) \cup Fix(f_2)$, this implies that $N(f) = 0$.   $\square$

## §8. MAPS ON SOME GEOMETRIC MANIFOLDS

In this section we will prove the main theorem for six types of geometric manifolds. The remaining two geometries are the $H^2 \times \mathcal{R}$ and $\widetilde{SL_2R}$, which will be treated in Section 9, together with manifolds admitting nontrivial JSJ splittings.

If $M$ is a manifold with toroidal boundary whose interior admits a complete hyperbolic structure, then we can identify $M$ with the $\epsilon$-thick part of a complete hyperbolic manifold. More explicitly, the interior $X$ of $M$ admits a complete hyperbolic structure with finite volume. Let $X_\epsilon$ be the set of all points $x$ in Int$X$ such that there is an embedded hyperbolic open ball of radius $\epsilon$ centered at $x$. There exists an $\epsilon_0 > 0$ such that for any $\epsilon \leq \epsilon_0$, $X_\epsilon$ is homeomorphic to $M$ and each boundary component of $X_\epsilon$ is a horosphere modulo a rank two abelian group of parabolic motions. (See [Th1, 5.11].) We identify $M$ with such an $X_\epsilon$.

**Lemma 8.1.** *Let $f$ be a map of a hyperbolic manifold $M$. Then $f$ is isotopic to a type 2 standard map $g$.*

*Proof.* Since $M$ is identified with $X_\epsilon$, the map $f : M \to M$ extends to a map $f' : X \to X$. By Mostow's Rigidity Theorem [Mo, p.54], there is a unique isometry $g' : X \to X$ which is homotopic to $f'$. Since $X$ retracts to $X_\epsilon = M$, the map $g = g'|_M$ is homotopic to $f$. By a recent result of Gabai, Meyerhoff and Thurston [Ga, GMT], $f$ is actually isotopic to the isometry $g : M \to M$.

If $g = id$, then by a small flow along a non vanishing vector field which is tangential on the boundary, we get a fixed point free map, and we are done. So assume that $g \neq id$.

Suppose each of $A_0, A_1$ is either a fixed point of $g$ or a boundary component of $M$, and $A_0$ is $g$-related to $A_1$ via the path $\gamma$. The universal covering $\widetilde{M}$ of $M$ is the hyperbolic 3-space $H^3$ with some horospheres removed. By lifting $\gamma$ to $\widetilde{M}$ we see that $\gamma$ is $(A_0, A_1)$ homotopic to a geodesic $\gamma'$, which is perpendicular to $A_i$ if it is a boundary component. Moreover, such geodesic $\gamma'$ is unique among all paths which are $(A_0, A_1)$ homotopic to $\gamma$. Since $g$ is an isometry, the path $g \circ \gamma'$ is also a geodesic, and

$$g \circ \gamma' \sim g \circ \gamma \sim \gamma \sim \gamma' \qquad \text{rel } (A_0, A_1).$$

So by the uniqueness of $\gamma'$, we must have $\gamma' = g \circ \gamma'$. In particular, $\gamma'$ is a path in $Fix(g)$. Therefore $g$ has FR-property. Since $g$ is a non identity orientation preserving isometry, $Fix(g)$ consists of mutually disjoint geodesics, so it is a properly embedded 1-manifold. A neighborhood of $Fix(g)$ is parametrized by a family of hyperbolic disk perpendicular to $Fix(g)$, each of which is preserved by the isometry $g$, hence $g$ preserved a normal structure of $Fix(g)$. It is known that isometries on hyperbolic manifolds are periodic [Th1], so in this case $g|_{\partial M}$ is a periodic map with isolated fixed points.   $\square$

**Corollary 8.2.** *If $M$ is a closed hyperbolic manifold, and $f : M \to M$ is an orientation preserving map, then $f$ is isotopic to a fixed point free map.*

*Proof.* By Lemma 8.1 we may assume that $f$ is a type 2 standard map, so $Fix(f)$ is a 1-manifold. Since $M$ is closed, $Fix(f)$ consists of circles, which, by Lemma 5.5, can be removed via an isotopy.   $\square$

**Lemma 8.3.** *Let $f$ be a map on the manifold $T^3$. Then $f$ is isotopic to a map $g$ with $\#Fix(g) = N(f)$.*

*Proof.* By Lemma 4.1 we may assume that $f$ is a linear map. Consider a lifting $\widetilde{f}$ of $f$ on the universal covering $\mathbb{R}^3$. Then $Fix(\widetilde{f})$ is an affine set, projecting to a component of $Fix(f)$, which must be compact. Since $Fix(\widetilde{f})$ is always connected, by Lemma 5.12 the components of $Fix(f)$ are mutually inequivalent.

Let $C$ be a component of $Fix(f)$. If dim $C = 3$, $f$ is an identity map. A small flow of $f$ along a nonvanishing vector field on $M$ (which exists because $\chi(M) = 0$), would eliminate all fixed points. If dim $C = 2$, $C$ is an affine torus or Klein bottle in $M$, cutting $M$ into pieces, each of which is a Euclidean manifold with boundary, so it must be an $I$-bundle over torus or Klein bottle. In any case there is a torus bundle structure of $M$ over 1-orbifold with $C$ a torus fiber. Since $f(C) = C$, $f$ is isotopic to a torus fiber preserving map which induces an orientation preserving map on the quotient space $S^1$, hence by Lemmas 7.3, $f$ is also isotopic to a fixed point free map.

Suppose dim $C = 1$. There is a lifting $\widetilde{f} : \mathbb{R}^3 \to \mathbb{R}^3$, which has a 1-dimensional fixed point set $\widetilde{C}$ covering $C$. Choose $\widetilde{C}$ as the $z$-axis of $\mathbb{R}^3$. The lifting $\widetilde{f}$ is then a linear map. After resizing we may assume that the tube $N(\widetilde{C})$ of radius 2 around the $z$-axis covers a tubular neighborhood $N(C)$ of $C$, with covering transformations generated by the length 2 translation along the $z$-axis. Consider the linear map $h = j \circ \widetilde{f}|_{\mathbb{R}^2} : \mathbb{R}^2 \to \mathbb{R}^2$, where $\mathbb{R}^2$ is the $xy$-coordinate plane of $\mathbb{R}^3$, and $j : \mathbb{R}^3 \to \mathbb{R}^2$ is the orthogonal projection. First notice that $h$ can not be an identity map, otherwise the line $\widetilde{f}(\mathbb{R}^2) \cap \mathbb{R}^2$ would be in $Fix(\widetilde{f})$, contradicting the assumption that $Fix(\widetilde{f}) = \widetilde{C}$. If $Fix(h) = \{0\}$, a small equivariant flow along the $z$ direction near $\widetilde{C}$ will eliminate $C$ from $Fix(f)$. So assume that $Fix(h)$ is 1-dimensional, say the $x$-axis. Assume that $h(0,1) = (a,b)$. Let $\rho_u$ be a positive rotation of $\mathbb{R}^2$ with angle $u$ if $a \leq 0$, and a negative rotation if $a > 0$. Let $r = (x,y)$, and let $||r||$ be the norm of $r$. Define an isotopy in the radius 1 tube around the $z$ axis by

$$\varphi_t(r,z) = (\rho_{t(1-||r||)}(r), z+t).$$

It induces an isotopy $\phi_t$ in $N(C)$ which can be extended over $M$ by identity. One can check that $g = \phi_1 \circ f$ has $Fix(g) = Fix(f) - C$.

We can thus isotop $f$ to a map $g$ whose fixed points are the isolated fixed points of $f$, and $ind(g,x) = ind(f,x)$ for all $x \in Fix(g)$. By Lemma 5.13 each isolated fixed point of $Fix(f)$ has nonzero index. Therefore, the number of fixed points of $g$ is exactly $N(f)$.  $\square$

**Lemma 8.4.** *If $f$ is an orientation preserving map on the manifold $M = M_{P(2,2)}$, then it is isotopic to a fixed point free map.*

*Proof.* By Theorem 4.5 $f$ is isotopic to an isometry with respect to some Euclidean metric of $M$. Since $f$ is orientation preserving, its fixed point set consists of disjoint circles, which can be removed by Lemma 5.5.  $\square$

**Proposition 8.5.** *If $M$ is a closed orientable 3-manifold which admits a Sol, Nil, $H^3$, $E^3$, $S^3$ or $S^2 \times S^1$ geometry, then any orientation preserving map $f$ on $M$ is isotopic to a map $g$ such that $\#Fix(g) = N(f)$.*

*Proof.* If $M$ admits a Sol or hyperbolic manifold, $T^3$, or $M_{P(2,2)}$, the result follows from Proposition 7.5, Lemmas 8.2, 8.3 and 8.4, respectively. So assume that $M$ is not such a manifold. For all the remaining manifolds, by Theorem 3.11 the map $f$ is isotopic to a fiber preserving map with respect to some Seifert fibration of $M$ which has spherical or Euclidean orbifold. If $X(M)$ is Euclidean, by Theorem 2.5 $\widehat{f}$ is isotopic to a periodic, reducible, or Anosov map. Note that this is trivially true if $X(M)$ is spherical because all maps on a spherical orbifold is isotopic to a periodic map.

If $\widehat{f}$ is Anosov, by Proposition 6.3 $f$ is isotopic to a Type 1 map, and by Lemma 5.7 $f$ can be further isotoped so that $N(f) = \#Fix(f)$, as desired. If $X(M)$ is spherical, by Proposition 6.4 $f$ is isotopic to a fixed point free map. If $\widehat{f}$ is periodic and $X(M)$ is Euclidean, by Proposition 6.3 $f$ is isotopic to a Type 2 standard map $g$. Since $M$ is closed, each component of $Fix(g)$ is a circle; hence by Lemma 5.5 we can isotop $g$ to a fixed point free map.

Now assume that $\widehat{f}$ is reducible, and let $C$ be a reducing curve on $X(M)$. Then $M$ also has the structure of a torus bundle over a 1-orbifold, with the torus $T$ over $C$ a torus fiber. Since $f(T) = T$, we can isotop $f$ to a torus fiber preserving map. The Seifert fibration provides an $S^1$ bundle structure on $T$, which is preserved by $f$, hence $f|_T$ is a reducible map. We can now apply Lemmas 7.3 and 7.4 to conclude that $f$ is isotopic to a map with no fixed points.  □

## §9. PROOF OF THE MAIN THEOREMS

**Theorem 9.1.** *Suppose $M$ is a closed orientable manifold which is either Haken or geometric, and $f : M \to M$ is an orientation preserving homeomorphism. Then $f$ can be isotoped to a homeomorphism $g$ with $\#Fix(g) = N(f)$.*

*Proof.* If $M$ admits a Sol, Nil, $E^3$, $S^3$ or $S^2 \times S^1$ geometry, the result follows from Proposition 8.5. Assume that $M$ is not such a manifold, and let $\mathcal{T}$ be the (possibly empty) set of tori of the Jaco-Shalen-Johannson decomposition. Then each component of $M - \text{Int}N(\mathcal{T})$ is either hyperbolic, or a twisted $I$-bundle over Klein bottle, or a Seifert fiber space with hyperbolic orbifold. Isotop $f$ so that it maps $N(\mathcal{T})$ homeomorphically to itself.

Suppose that $P$ is a Seifert fibered component of $M - \text{Int}N(\mathcal{T})$ such that $f(P) = P$. By [Sc, Theorem 3.9], $f|_P$ is isotopic to a fiber preserving map. By Lemma 1.10, we can find a set of vertical tori $\mathcal{T}_P$ in $P$, cutting $P$ into pieces with hyperbolic or $D(2,2)$ orbifold, and a fiber preserving isotopy of $f$, so that after isotopy, the restriction of $f$ on each invariant piece has periodic or pseudo-Anosov orbifold map. Adding all such $\mathcal{T}_P$ to $\mathcal{T}$, we get a collection of tori $\mathcal{T}'$, such that

(i) $f(N(\mathcal{T}')) = N(\mathcal{T}')$;

(ii) Each component $M_i$ of $M - \text{Int}N(\mathcal{T}')$ is either hyperbolic, or a Seifert fiber space such that $X(M_i)$ is a hyperbolic orbifold or $D(2,2)$.

(iii) If $f$ maps a Seifert fibered component $M_i$ to itself, and if $M_i$ has hyperbolic orbifold, then $f$ is fiber preserving, and the orbifold map $\widehat{f}$ on $X(M_i)$ is either periodic or pseudo-Anosov.

We can now apply Proposition 4.3, Lemma 5.1 and Lemma 6.1 to each component $M_i$ of $M - \text{Int} N(\mathcal{T}')$ which is mapped to itself by $f$, isotoping $f$ so that $f|_{M_i}$ is a standard map for all $M_i$. After that, use Lemma 6.2 to further isotop $f$, rel $\partial N_j$, on each component $N_j$ of $N(\mathcal{T}')$, so that it is a standard map on $N_j$.

By the definition of standard maps, $Fix(f)$ intersects each of $M_i$ and $N_j$ in points and 1-manifolds, so $Fix(f)$ is a disjoint union of points, arcs, and circles. Moreover, at each isolated point of $Fix(f)$ the map $f$ is of pseudo-Anosov-flip type, and at each arc or circle it preserves some normal structure.

We want to show that different components of $Fix(f)$ are not equivalent. If this is not true, there is a path $\alpha$ connecting different components $C_0, C_1$ of $Fix(f)$, such that $f \circ \alpha \sim \alpha$ rel $\partial$. Denote by $T''$ the set of tori $\partial N(\mathcal{T}')$. Among all such $\alpha$, choose one such that $|\alpha - \mathcal{T}''|$ is minimal. In below we will find another such curve $\alpha''$ with $|\alpha'' - \mathcal{T}''| < |\alpha - \mathcal{T}''|$, which would contradict the choice of $\alpha$.

Let $D$ be a disk, and let $h : D \to M$ be the homotopy $\alpha \sim f \circ \alpha$ rel $\partial$. We may assume that $h$ is transverse to $\mathcal{T}''$, and $|h^{-1}(\mathcal{T}'')|$ is minimal among all such $h$. Then $h^{-1}(\mathcal{T}'')$ consists of a properly embedded 1-manifold on $D$, together with possibly one or two isolated points mapped to the ends of $\alpha$. Clearly, $\mathcal{T}''$ is $\pi_1$-injective in $M$, so one can modify $h$ to remove all circles in $h^{-1}(\mathcal{T}'')$. Note that $h^{-1}(\mathcal{T}'')$ must contain some arcs, otherwise $\alpha$ would lie in some $M_i$ or $N(T_j)$, which is impossible because the restriction of $f$ in each piece has FR-property.

Now consider an outermost arc $b$ in $h^{-1}(\mathcal{T}'')$. Let $\beta = h(b)$. The ends of $\beta$ can not both be on $\alpha$, otherwise we can use the outermost disk to homotope $\alpha$ and reduce $|\alpha - \mathcal{T}''|$, contradicting the choice of $\alpha$. Since $f$ is a homeomorphism, the same thing is true for $f \circ \alpha$. Therefore, $\beta$ has one end on each of $\alpha$ and $f \circ \alpha$.

The arc $b$ cuts off a disk $\Delta$ on $D$ whose interior is disjoint from $h^{-1}(\mathcal{T}'')$. The boundary of $\Delta$ gives rise to a loop

$$h(\partial \Delta) = \alpha_1 \cup \beta \cup (f \circ \alpha_1)^{-1},$$

where $\alpha_1$ is a subpath of $\alpha$ starting from an end point $x$ of $\alpha$. Let $T$ be the torus in $\mathcal{T}''$ which contains $\beta$. Then the restriction of $h$ on $\Delta$ gives a homotopy

$$\alpha_1 \sim f \circ \alpha_1 \quad \text{rel } (x, T).$$

Since $f$ has FR-property on each component of $M - \mathrm{Int} N(\mathcal{T})$ and $N(\mathcal{T})$, by definition there is a path $\gamma$ in $Fix(f)$ such that $\gamma \sim \alpha_1$ rel $(x, T)$. Since $\gamma$ is in $Fix(f)$, the path $\alpha' = \gamma^{-1} \cdot \alpha$ has the property that

$$f \circ \alpha' = (f \circ \gamma^{-1}) \cdot (f \circ \alpha) = \gamma^{-1} \cdot (f \circ \alpha) \sim \gamma^{-1} \cdot \alpha = \alpha' \quad \text{rel } \partial.$$

Since $\alpha_1 \sim \gamma$ rel $(x, T)$, the path $\gamma^{-1} \cdot \alpha_1$ is rel $\partial$ homotopic to a path $\delta$ on $T$. Write $\alpha = \alpha_1 \cdot \alpha_2$. Then

$$\alpha' = \gamma^{-1} \cdot \alpha = (\gamma^{-1} \cdot \alpha_1) \cdot \alpha_2 \sim \delta \cdot \alpha_2 \quad \text{rel } \partial.$$

By a small perturbation on $\delta$, we get a path $\alpha'' \sim \alpha'$ rel $\partial$ such that $|\alpha'' - \mathcal{T}''| < |\alpha - \mathcal{T}''|$. Since

$$f \circ \alpha'' \sim f \circ \alpha' \sim \alpha' \sim \alpha'' \quad \text{rel } \partial,$$

this contradicts the minimality of $|\alpha - \mathcal{T}''|$.

So far we have proved that each component of $Fix(f)$ is a Nielsen class of $f$. Our next goal is to show that we can isotop $f$ in a neighborhood of each arc component $C$ of $Fix(f)$ so that the class $C$ becomes a single point, or is completely deleted if $ind(f, C) = 0$.

Notice that $C$ is the union of several arcs $C_k$, each being a component of some $Fix(f|_{M_i})$ or $Fix(f|_{N_j})$, and the ends of $C$ are isolated points of, say, $Fix(f|_{M_0})$ and $Fix(f|_{M_1})$. Since $f$ is a standard map on each piece, we can identify $N(C)$ with a subset of $\mathbb{R}^3$ as follows:

(1) $N(C) = [-3, 3] \times D$, where $D = [-2, 2] \times [-2, 2]$.

(2) $C = [-1, 1] \times (0, 0)$.

(3) For $t \in [-1, 1]$, $f$ maps $t \times D_1$ into $t \times D$, where $D_1 = [-1, 1] \times [-1, 1]$ is a smaller neighborhood of $(0, 0)$. This is because $f$ preserves a normal structure on each arc component of $Fix(f|_{M_i})$ or $Fix(f|_{N_j})$.

(4) On $[-2, -1] \times D$, the map has the property that $f(x, y, z) = (\varphi(x, y), -z)$, where $\varphi(x, y)$ has the property that it switches the two sides of the $xy$-plane on the two sides of the $x$-axis, i.e, if $\varphi(x, y) = (x', y')$ then $y' > 0$ if and only if $y < 0$. This follows from the fact that the point $(-1, 0, 0)$ is an isolated point of $Fix(f|_{M_0})$, hence is of pseudo-Anosov-flip type. Similarly for $[1, 2] \times D$.

Consider the ball $B_1 = [-2, 2] \times D_1$. Define a map $h : B_1 \to \mathbb{R}^3$ by $h(v) = v - f(v)$. It follows from the above discussion that $f|_{B_1}$ have the following properties.

(5) $h(v)$ is parallel to the $x$ axis if and only $v = (x, 0, 0)$ for some $x$. Moreover, for all $-2 \leq x < -1$, $h(v)$ points to the same direction of the $x$-axis. Similarly for all $1 < x \leq 2$.

We can now compute $ind(f, C)$. Put $e = (1, 0, 0)$. By (5) there are only two points on $\partial B_1$ that are mapped to $\pm e$ by the map $\bar{h}(v) = h(v)/||h(v)||$, namely the points $\pm 2e$. It follows that the degree of the map $\bar{h} : \partial B_1 \to S^2$ is either $0$ or $\pm 1$, depending on whether $\bar{h}(2e) = \bar{h}(-2e)$ or not.

If $\bar{h}(2e) = \bar{h}(-2e) = e$, say, then $v - f(v)$ is never a positive multiple of $-e$. Define an isotopy supported in a small neighborhood of $C$ by flowing along the direction of $-e$. This isotopy will modify $f$ to a map $g$ with no fixed point in $B_1$.

If $\bar{h}(2e) = e$ and $\bar{h}(-2e) = -e$, say, we can define an isotopy supported in a small neighborhood of $C$ by flowing each point $(x, y, z)$ towards $(0, y, z)$. This will isotop $f$ to a new map $g$ with the origin as the only fixed point in $B_1$. Since $ind(f, C) \neq 0$, the origin is an essential fixed point of $g$. This completes the proof that we can isotop $f$ so that each arc component $C$ of $Fix(f)$ becomes a single point, or totally disappears, depending on whether $C$ is an essential fixed point class of $f$.

If $x$ is an inessential fixed point of $f$, we can remove it using Lemma 5.7, and if $C$ is a circle component of $Fix(f)$, we can remove it using Lemma 5.5. This, together with the above isotopies near each arc component of $Fix(f)$, will deform $f$ to a new map $g$ such that each component of $Fix(g)$ is an essential fixed point class. Therefore, $\#Fix(g) = N(f)$.  □

Define an orbifold $X(M)$ to be *small* if it is a sphere with a total of at most three holes or cone points, or a projective plane with a total of 2 holes or cone points. Note that any map on a small orbifold is isotopic to a periodic map.

**Theorem 9.2.** *Let $M$ be a closed orientable 3-manifold. Then any orientation preserving homeomorphism $f$ on $M$ is isotopic to a fixed point free homeomorphism, unless some component of the JSJ decomposition of $M$ is a Seifert fiber space with big orbifold.*

*Proof.* This follows from Proposition 7.5 if $M$ is a Sol manifold, and from Corollary 8.2 if $M$ is hyperbolic. If $M$ is a Seifert fiber space with $X(M)$ a sphere with at most 3 cone points, then by Theorem 3.11, $f$ is isotopic to a fiber preserving map. It follows immediately from Proposition 6.4 that $f$ is isotopic to a fixed point free map.

Now suppose $M$ has a nonempty JSJ decomposition. Let $\mathcal{T}$ be the decomposition tori, and isotop $f$ so that $f(N(\mathcal{T})) = N(\mathcal{T})$. By the same argument as above, we may assume that the restriction of $f$ on each component $M_i$ of $M - \text{Int}N(\mathcal{T})$ is a type 2 standard map.

Let $T$ be a component of $\mathcal{T}$, and assume that $f$ maps $N(T)$ to itself. Suppose $f$ exchange the two boundary components of $N(T)$. Then $f^2$ maps each piece $M_i$ adjacent to $N(T)$ to itself. Since $M_i$ is either hyperbolic or Seifert fibered, $f^2|_{M_i}$ is isotopic to a periodic or Seifert fiber preserving map, so the restriction of $f^2$ to the boundary of $N(T)$ can not be isotopic to an Anosov map. By Lemma 7.2(1), $f|_{N(T)}$ is rel $\partial N(T)$ isotopic to a fixed point free map, which is automatically a type 2 map. If $f$ maps each boundary component of $N(T)$ to itself, then $f$ is already a type 2 standard map on the two pieces adjacent to $N(T)$, so by Lemma 7.2(2) $f|_{N(T)}$ is also rel $\partial N(T)$ isotopic to a type 2 map.

After the above isotopy, we get a map $f$ which is of type 2 on each $M_i$ and $N(T)$. Thus each component of $Fix(f)$ is a circle, which can be removed by Lemma 5.5.    $\square$

DEPARTMENT OF MATHEMATICS, PEKING UNIVERSITY, BEIJING 100871, CHINA
*E-mail address*: jiangbj@sxx0.math.pku.edu.cn

DEPARTMENT OF MATHEMATICS, PEKING UNIVERSITY, BEIJING 100871, CHINA; and MSRI, 1000 CENTENNIAL DRIVE, BERKELEY, CA 94720-5070
*E-mail address*: wangsc@msri.org

DEPARTMENT OF MATHEMATICS, UNIVERSITY OF IOWA, IOWA CITY, IA 52242; and MSRI, 1000 CENTENNIAL DRIVE, BERKELEY, CA 94720-5070
*E-mail address*: wu@math.uiowa.edu, wu@msri.org